\documentclass[11pt]{amsart}
\usepackage{amsmath,amsthm, amscd, amssymb, amsfonts, xypic, xcolor}
\usepackage{graphicx}
\usepackage{tikz-cd}
\theoremstyle{plain}

\newtheorem{theorem}{Theorem}[section]
\newtheorem{lemma}[theorem]{Lemma}

\newtheorem{proposition}[theorem]{Proposition}

\theoremstyle{definition}
\newtheorem{definition}[theorem]{Definition}
\newtheorem{corollary}[theorem]{Corollary}
\newtheorem{example}[theorem]{Example}

\newcounter{remarks}
\theoremstyle{remark}
\newtheorem{remark}[remarks]{Remark}

\newcommand{\codim}{\textrm{codim}}
\newcommand{\rank}{\operatorname{rank}}

\newcommand{\defect}{{\operatorname {def}}}

\newcommand{\baseRing}[1]{\ensuremath{\mathbb{#1}}}
\newcommand{\Z}{\baseRing{Z}}

\newcommand{\C}{\baseRing{C}}

\newcommand{\Q}{\baseRing{Q}}

\newcommand{\CP}{\baseRing{P}}

\newcommand{\Script}[1]{\ensuremath{{\mathcal{#1}}}}
\newcommand{\ScS}{\Script{S}}

\newcommand{\bp}{\begin{proof}}
\newcommand{\ep}{\end{proof}}
\newcommand{\beq}{\begin{equation}}
\newcommand{\eeq}{\end{equation}}

\newcommand{\beqs}{\begin{equation*}}
\newcommand{\eeqs}{\end{equation*}}

\newcommand{\beas}{\begin{eqnarray*}}
\newcommand{\eeas}{\end{eqnarray*}}

\newcommand{\btheorem}{\begin{theorem}}
\newcommand{\etheorem}{\end{theorem}}

\newcommand{\bl}{\begin{lemma}}
\newcommand{\el}{\end{lemma}}

\newcommand{\benum}{\begin{enumerate}}
\newcommand{\eenum}{\end{enumerate}}

\renewcommand{\k}{{\mathbb K}}

\newcommand{\fs}{\mathfrak s}

\usepackage{todonotes}

\newcommand{\csig}{\ensuremath{\underline{\mathbf \sigma}}}
\newcommand{\cset}{\underline{S}}
\newcommand{\cdiff}{\underline{D}}

\numberwithin{equation}{section}
\begin{document}

\title[Flags, Circuits, Matrices and Cayley Configurations]{Non-splitting flags, Iterated Circuits, $\csig$-matrices and
Cayley configurations}

\author[E. Cattani ]{Eduardo Cattani}
\address{Department of Mathematics and Statistics, 
University of Massachusetts Amherst,
Amherst MA 01003-9305,
USA }
\email{cattani@math.umass.edu}

\author[A. Dickenstein]{Alicia Dickenstein}
\address{Dto.\ de Matem\'atica, FCEN, Universidad de Buenos Aires, and IMAS (UBA-CONICET), Ciudad Universitaria, Pab.\ I, 
C1428EGA Buenos Aires, Argentina}
\email{alidick@dm.uba.ar}
\urladdr{http://mate.dm.uba.ar/~alidick}

\thanks{AD was partially supported by UBACYT 20020100100242, CONICET PIP 20110100580, and ANPCyT PICT 2013-1110, Argentina}  

\date{}

\maketitle

\begin{center}  {\em Dedicated to Bernd Sturmfels
on the occasion  of his $60$th birthday}
\end{center} 
\bigskip

\begin{abstract}
We explore four approaches to the question of defectivity for a complex projective toric variety $X_A$ associated with an integral configuration $A$.  
The explicit tropicalization of the dual variety $X_A^\vee$ due to Dickenstein, Feichtner, and Sturmfels allows for the computation of the defect in terms of 
an affine combinatorial invariant $\rho(A)$.  We express $\rho(A)$ in terms of affine invariants $\iota(A)$ associated to Esterov's iterated circuits and $\lambda(A)$, 
an invariant defined by Curran and Cattani in terms of a Gale dual of $A$. Thus we obtain formulae for the dual defect in terms of iterated circuits and Gale duals. 
An alternative expression for the dual defect of $X_A$ is given by Furukawa-Ito in terms of Cayley decompositions of $A$. 
We give a Gale dual interpretation of these decompositions and apply it to the study of defective configurations.

\end{abstract}

\section{Introduction}\label{sec:intro}

Given a complex projective variety $X$, its dual $X^\vee$ is defined as the closure in the dual 
projective space of all hyperplanes tangent to $X$ at a smooth point.  
It is classically known that generically $X^\vee$ is a hypersurface (\cite{GKZ}, Corollary~1.2).  
If $\codim\,X^\vee >1$,  $X$ is said to be {\em defective} and the quantity 
$$
\defect(X) \ :=\ \codim\,X^\vee -1
$$
is called the {\em dual defect} of $X$.  If $X$ is irreducible and non-defective, then $X^\vee$ is irreducible and the polynomial 
defining $X^\vee$ is known as the {\em discriminant} of $X$.

Suppose now that $A = \{a_1,\dots,a_n\} \subset \Z^e$ is a configuration (of not necessarily distinct points) which is 
{\em homogenous} in the sense that the $n$-tuple $(1,\dots,1)$ lies in 
the row span of $A$ viewed as the $e\times n$ matrix whose columns are the $a_i$'s.    Then $A$ defines a projective 
toric variety $X_A \subset \CP^{n-1}$~\cite{GKZ} rationally parametrized by monomials with exponents in $A$.
The dimension $\dim(X_A)$ equals $\rank(A)-1$, which is equal to the affine dimension $d(A)$ of $A$. 

We will also assume that $A$ 
is not a {\em pyramid}; i.e. no affine hyperplane contains all points of 
$A$ except one or, equivalently, the dual variety $X_A^\vee$ is not contained in a hyperplane. We discuss at the end of 
Section~\ref{sec:right=} the extension of our definitions and results to the pyramidal case.

The dual defect of $X_A$ is an affine invariant of the configuration $A$ which has been studied from various points of view.  
In this paper we will discuss four such approaches to the  dual defect of projective toric varieties. 
Let us briefly describe them postponing for a moment the detailed definitions:
 
\begin{enumerate}
\item[i)] 
In \cite{DFS} Dickenstein, Feichtner, and Sturmfels described the tropicalization of the dual variety $X_A^\vee$. 
This leads to the computation of the dual defect of $X_A$, explicitly, 
in terms of certain $n(A)\times n(A)$ matrices $M_{\csig}(A)$, where $n=n(A)$ is the cardinality of $A$. 
These matrices are constructed from the matrix $A$ and a maximal chain $\csig$ 
of support of vectors in the kernel of $A$ (see Definition~\ref{def:DFS}).  Indeed \cite[Corollary~4.5]{DFS} gives the equality 
\beq\label{eq:defect_rho}
\defect(X_A) = n(A) - 1 - \rho(A),
\eeq
 where 
$$
 \rho(A) = \max_{\csig} {\rm rank} (M_{\csig}(A))),
$$
We refer to $M_{\csig}(A)$ as a ${\csig}$-matrix.  We point out that a key ingredient in the computation of the dual 
defect is the Horn-Kapranov parametrization map, also used by Forsg{\aa}rd~\cite{Fo}, who applies it to real configurations.

\medskip

\item[ii)] Recall that a point configuration $Z$ is called a \emph{circuit}  if it is minimally affinely dependent 
i.e. it is affinely dependent but every proper subset is affinely independent. It is not hard to show that if $A$ contains a circuit $Z$ with $d(Z)=d(A)$ then $X_A$ is non-defective.  
In \cite{E2}, Esterov shows that if $\defect(X_A) =0$ then $A$ might not contain a circuit of  affine dimension $d(A)$, 
but it will necessarily contain what he calls an {\em iterated circuit} 
with this maximal possible affine dimension. We generalize this notion in Definition~\ref{def:itercircuit}.

\medskip

\item[iii)] We can associate to the configuration $A$ a Gale dual configuration $B$.
Let $B \in \Q^{n \times m}$, where $m=m(A) = \dim \ker_\Q(A)$, be a matrix whose columns 
are a basis of $\ker_\Q A \subset \Q^n$. The rows of $B$ define a linear 
matroid which we also denote by $B$ and which gives a choice of Gale dual configuration (see Section~\ref{sec:right=}). 
If $A$ is homogeneous, then $B$ satisfies a dual homogeneity condition:
$$\sum_{b\in B} b = 0.$$
Associated to the matroid $B$ is its lattice of flats  (see Section~\ref{sec:right=}).  We say that a chain of flats
$$F_1 \subset F_2 \subset \cdots \subset F_\ell\ ;\quad \dim(F_j) = j,$$ 
is {\em non-splitting} if $\sum_{b\in F_1} b \not=0$ and 
$\sum_{b\in F_{j+1}} b \not\in L(F_j)$, $j\geq 1$,
where $L(F_i)$ denotes the  linear span of the vectors in $F_i$.  We call $\ell$ the length of the flag and 
denote by $\lambda(B)= \lambda(A)$ the maximal length of a  non-splitting flag of flats in $B$.  
It is shown in \cite{CC} that $X_A$ is defective if and only if $\lambda(B) < m -1$.   
This generalizes the work of Dickenstein and Sturmfels \cite{DS} in the case $m=2$.

\medskip

\item[iv)] Every defective configuration admits a Cayley decomposition (see Definition~\ref{def:cayley}) 
but this is not enough to characterize defectivity.  However, Furukawa and Ito \cite{FI} 
show that among all possible Cayley decompositions of $A$ there is a distinguished one, 
which we will call the {\em FI-decomposition}, and it is then possible
 to compute $\defect(X_A)$ in terms of simple invariants of this decomposition (see Section~\ref{sec:cayley}.)  
 Their approach is fundamentally different to the tropical approach used in \cite{DFS}.  
 
\end{enumerate}

The purpose of this paper is twofold.  First of all, we extend the computation of $\defect(X_A)$
 to the settings described in ii) and iii).  In Section~\ref{sec:left=}, 
we extend Esterov's notion of iterated circuit to the defective case and define an invariant $\iota(A)$ 
which is the maximal rank of an iterated circuit $I \subset A$. We show that
\beq\label{eq:itercirc}
n(A) - 1 - \rho(A) = d(A) - \iota(A).
\eeq
This allows us to compute the defect of $X_A$ in the context of Esterov's work. 
On the other hand, in Section~\ref{sec:right=}
we prove  that 
\beq\label{eq:flags}n(A) - 1 - \rho(A) = m(A) -1 - \lambda(A),
\eeq
and, consequently, the right-hand side computes the dual defect of $X_A$ in terms of the Gale dual $B$.

We should point out that, though inspired by the geometric case of integral configurations and the 
associated toric varieties over $\mathbb C$, the proof of these 
results is linear-algebraic in nature and  the above equalities hold for finite configurations 
$A = \{a_1,\dots,a_n\} \subset \k^e,$ defined over an arbitrary field $\k$ of characteristic zero. Thus, we will work in this generality.

\medskip

We give all the definitions about iterated circuits and $\csig$-matrices in Section~\ref{sec:left=}, where 
we prove equality~(\ref{eq:itercirc}) in Theorem~\ref{th:itercirc}. In Section~\ref{sec:right=} we introduce necessary
concepts about Gale duality  and homogenizations and we prove equality~(\ref{eq:flags}) in Theorem~\ref{th:gale}.
The expression (\ref{eq:defect_rho}) for $\defect(X_A)$ given in \cite{DFS} is valid provided $A$ is not a pyramid.  At the end of Section~\ref{sec:right=}
 we include a brief discussion of how our results can be applied in the pyramidal case.

Finally,  Section~\ref{sec:cayley} should be considered as the second part of this paper.  We begin by describing the results of Furukawa and Ito in \cite{FI}
 which allow us to compute the dual defect of complex toric varieties $X_A$ in terms of the family of Cayley decompositions of $A$.  We then translate this to 
 the Gale dual setting and use this formulation to obtain results in low codimension and for configurations with ``large" defect.

 \medskip

\noindent{\bf Acknowledgment:} We are indebted to Atsushi Ito for fruitful discussions and for allowing us to quote
his counterexample~\ref{ex:char2}.  We thank the referees for their detailed and useful comments to improve
our manuscript.   Almost 30 years ago, at an NSF Geometry Institute held at Amherst College, Bernd Sturmfels introduced 
us to many of the topics discussed here. We are grateful for his inspiration and his friendship throughout this time.

\bigskip

\section{Iterated Circuits and $\csig$-Matrices: Proof of (\ref{eq:itercirc})}\label{sec:left=}

In this section we generalize Esterov's notion of {\em iterated circuits} (see Definition~\ref{def:itercircuit}) to include
the case of defective configurations and prove identity (\ref{eq:itercirc}) which allows us to compute the dual defect of a toric variety in terms of iterated circuits. 
This is the content of   the main result in this section, Theorem~\ref{th:itercirc}, which follows from Theorem~\ref{prop:bij} 
where we make explicit the relationship between iterated circuits and  chains of
supports.

We begin by establishing the notation that we will use throughout the paper.  We will also give careful definitions of the objects mentioned in the Introduction.  
As has been noted already, we will work with a finite configuration $A = \{a_1,\dots,a_n\} \subset \k^e,$ defined over an arbitrary field $\k$ of characteristic zero. 

There are three basic affine quantities associated with $A$, namely, its cardinality $n(A)$, its affine dimension $d(A)$ 
over $\k$ (that is, the $\k$-dimension of its {\em affine span} $L_{\rm aff}(A))$, and 
the dimension $m(A)= n(A)-(d(A)+1)$ of the space $R_{\rm aff}(A)$ of $\k$-affine 
relations among the points in $A$. We also refer to $m(A)$ as the {\em codimension} of $A$.   If there is no possibility of confusion, we will simply write $n$, $d$, 
and $m$ for the quantities $n(A)$, $d(A)$ and $m(A)$.

As before, we may identify $A$ with the $e \times n$ matrix whose columns are the $a_i$'s. We will say that $A$ is homogeneous 
if the $n$-tuple $(1,\dots,1)$ lies in the row span of $A$. Up to a linear transformation, this is equivalent to saying that $A$ 
is contained in an affine hyperplane not passing through the origin. 
If $A$ is homogeneous then 
$R_{\rm aff}(A) = \ker_\k(A)$.

In what follows we also need to consider non-homogeneous configurations.  If $A$ is not homogeneous we associate to it two homogeneous configurations:
in case $A$ is not homogeneous, we denote by $\bar{A}$ the homogeneous configuration
\begin{equation}\label{eq:overline}
\bar{A} \, = \{ (1,a_1), \dots, (1, a_{n})\} \subset \k^{e+1},
\end{equation}
which is affinely equivalent to $A$.   We also set 
\begin{equation}\label{def:Chom}
A^h  \, = \, \{(1, 0, \dots, 0)\} \cup \bar{A} \, = \, \overline{(0,\dots,0) \cup A}\subset \k^{e+1}.
\end{equation}
We will refer to $A^h$ as the {\em homogenization} of $A$. When $A$ is homogeneous, it is enough to take $\bar{A}=A$.

\begin{remark} \label{rem:circuit} We note that if $Z$ is a circuit, that is a minimally affinely dependent set of points, then $Z$ is not a pyramid and 
$n(Z) = d(Z)+2$ (or, equivalently, $m(Z)=1$).  In particular, a zero-dimensional circuit is a configuration of cardinality two with a repeating point.  
Observe also that if $W \subset \k^e$ is minimally {\em linearly} dependent over $\k$, then it is a circuit provided it is 
homogeneous and, otherwise, its homogenization $W^h$ is a circuit.
\end{remark}

We will denote  the {\em $\k$-linear span} of a subset $S$ in a $\k$-vector space by  $L(S)$.

\smallskip

The following is a generalization of the definition of iterated circuit given
	by Esterov~\cite{E1,E2}.

\begin{definition}\label{def:itercircuit}
Let $A$ be a non-pyramidal homogeneous point configuration over $\k$ of maximal rank.
An {\em iterated circuit of length $p$} in $A$ is a subconfiguration $I \subset A$
together with a {\em partition} into nonempty subconfigurations
$$ I =  I_1 \cup \cdots \cup I_p,$$
such that:
\begin{enumerate}
\item[i)] $I_1$ is a circuit;
\item[ii)] For every $j$, $1\leq j\leq p-1$, write 
$L_{j}= L(I_1 \cup \cdots \cup I_{j})$,   and consider the natural projection
$\pi_{j} \colon L(I) \to L(I)/ L_j$. 
Then, either $\pi_{j}(I_{j+1})$  is a circuit with affine dimension $0$, or $\pi_{j}$ acts injectively on $I_{j+1}$.
Moreover, in the latter case,  the configuration $\pi_{j}(I_{j+1})$ is minimally linearly dependent (cf. Remark~\ref{rem:circuit}).
\end{enumerate}
\end{definition}
 
Given an iterated circuit $I =I_1 \cup \cdots \cup I_p \subset A$, let
\begin{equation}\label{eq:eta}
\eta(I) \, = \,  \sum_{j=1}^p \,d(I_j).
\end{equation}
Note that $\eta(I)$ equals $d(I)$  minus the number of indices
$j=1, \dots p-1$ for which the projection $\pi_{j}(I_{j+1})$ is homogeneous.

We now define
\begin{equation}\label{eq:iota}
\iota(A) \, = \,  {\rm max} \{\eta(I), \, I \subset A \hbox{ \ iterated circuit} \}.
\end{equation}

 We point out that in \cite{E2} Esterov reserves the name iterated circuit for iterated circuits $I$ with $\eta(I)=d(A)$. 
 In this case $\eta(I) = d(I)$ and none of the configurations   $\pi_{j}(I_{j+1})$, 
  $j=1, \dots p-1$, may be homogeneous.  For a non-homogeneous configuration $A$ we set   $\iota(A) = \iota(\bar{A})$.
  
  \begin{example}\label{ex:itcirc}
We present here two simple examples to illustrate Definition~\ref{def:itercircuit}. The figures below depict 
 non-homogeneous configurations $A'$ and we take $A = \overline{A'}$

\begin{figure}[h]
\label{fig:iter}
\includegraphics[width=5cm]{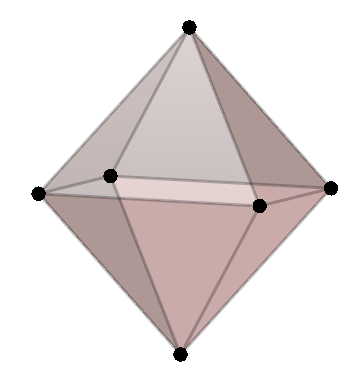} \qquad
\includegraphics[width=5cm]{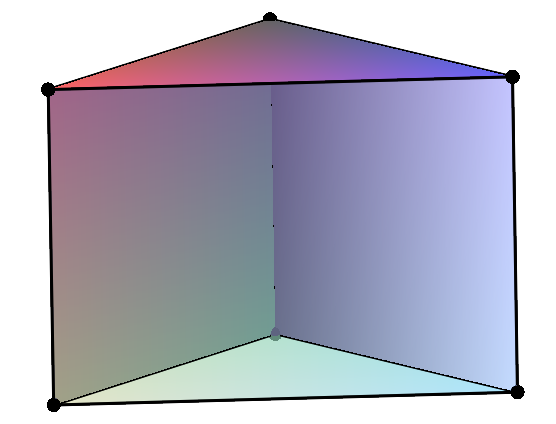}
\caption{Octahedron (left) and prism (right) configurations}
\end{figure}

In the octahedron configuration, we have 
\beq\label{octahedron}
A' = \{ (1,0,0),  (-1,0,0), (0,1,0), (0,-1,0), (0,0,1),(0,0,-1)\}.
\eeq 
It is clear that the four points in each of the coordinate planes define two-dimensional circuits. Setting 
 $I_1=\{ (1,1,0,0),  (1,-1,0,0), (1,0,1,0), (1,0,-1,0)\}$ and $I_2=\{ (1,0,0,1),(1,0,0,-1)\}$ 
 we obtain a maximal dimension iterated circuit in $A$.
Indeed, $I_1$ is a homogeneous circuit of dimension $2$ and $\pi_1(I_2)$ consists of two different points in dimension $1$,
 which is a minimally linearly dependent configuration and, consequently,  its homogenization is a circuit as noted in 
 Remark~\ref{rem:circuit}. Clearly, $\eta(I) = 2 + 1 =3$ and $A$ is non-defective.
 
Consider next the prism:
 \beq\label{prism} A' = \{(1,0,0), (1,0,1), (0,1,0), (0,1,1),  (0,0,0), (0,0,1)\}.
 \eeq
 Once again, there are three, two-dimensional circuits, but none of them may be extended to a 
 three-dimensional iterated circuit in $A$. For example, if we start with the circuit
$I_1 = \{ (1,0,1,0), (1,0,1,1),  (1,0,0,0), (1,0,0,1))\}$ then $\pi_1(A \setminus I_1)$ consists of two equal points and thus, is a minimally linearly 
dependent configuration which is a homogeneous zero-dimensional circuit.  
Thus $I = I_1 \cup (A \setminus I_1)$ is an iterated circuit but $\eta(I) = 2$.  Indeed, $\iota(A) =2$, and 
the configuration $A$ is defective. 
  
  \end{example}
  
\smallskip

We now recall the invariant introduced  by Dickenstein, Feichtner, and Sturmfels in~\cite{DFS}, which
 in the case of integral configurations yields the dual defect of the associated complex toric variety.
\begin{definition} \label{def:DFS}
Given a homogeneous non-pyramidal configuration $A$ over $\k$,
we denote by $\ScS(A)\subset \{0,1\}^{n(A)}$  the geometric lattice whose elements are the   supports, 
  ordered by inclusion, of the vectors in $R_{\rm aff} (A) = {\rm ker}_\k({A})$.  
  Given a maximal chain of  supports
  $$\csig = \{\sigma_1 \prec \sigma_2 \prec \cdots \prec \sigma_{m(A)}\} $$
  in $\ScS(A)$ we define the $n(A)\times n(A)$  matrix
  $$M_{\csig}(A) = \left(\begin{array}{c}  A \\\sigma_1 \\\vdots \\ \sigma_{m(A)}\end{array}\right).$$
 We will refer to $M_{\csig}(A)$ as a $\csig$-matrix associated with $A$ and set:
 \begin{equation}\label{eq:rho}
 \rho(A) = \max_{\csig} {\rm rank} (M_{\csig}(A))),
 \end{equation}
 where $\csig$ runs over all maximal chains of  elements in $\ScS(A)$. 
 
 For a non-homogeneous configuration $A$ we set $\rho(A) = \rho(\bar{A})$.
 \end{definition}

The implicit assertion in Definition~\ref{def:DFS} that the length of a maximal chain of supports is $m(A)$ 
follows from the fact that $\k$ is infinite, since we assume that ${\mathrm{char}}(\k) = 0$. 
  Hence, given any two vectors $v_1, v_2\in \k^e$, there exists
a $\k$ linear combination $v=\lambda_1 v_1 + \lambda_2 v_2$ such that the support of $v$ is the union of the supports of $v_1$ and $v_2$.  
Thus, the lenght of a maximal chain of supports in  ${\rm ker}_\k({A})$ equals $\dim ({\rm ker}_\k({A})) = m(A)$. 

\medskip

Before stating Theorem~\ref{th:itercirc},
 we prove some preliminary results on iterated circuits and $\csig$-matrices.

Given an iterated circuit 
$I=I_1 \cup \dots \cup I_p,$
an iterated circuit $I' = I'_1 \cup \dots \cup I'_{p'}$ is called an {\em extension} of 
$I$ if  $p \le p'$ and $I'_j = I_j$, for $j=1,\dots, p$.  Clearly, if $I'$ is an extension of $I$ then $\eta(I') \geq \eta(I)$.
We have:

\begin{lemma} \label{lem:imax}
Let $A$ be a homogeneous non-pyramidal configuration over $\k$. Then 
 any iterated circuit $I=I_1 \cup \dots \cup I_p$ has an extension $I'$ such that $d(I')=d(A)$.  
\end{lemma}

\begin{proof} Let $I'$ be an extension of $I$, maximal with respect to affine dimension.  
If $L(I') \not= L(A)$,  consider the natural projection
$\pi \colon L(A) \to L(A)/L(I')$. Then $\pi(A \setminus L(I'))$ 
 or its homogenization contain no circuits.  This means that $\pi(A \setminus L(I'))$ consists
  of linearly independent vectors, contradicting the fact that $A$ is non-pyramidal.  
\end{proof}

We now introduce some notation which will help us to better understand
the relationship between iterated circuits and  maximal chains of supports in $\ker_\k(A)$.  

Denote by $[n]$ the set $\{1,\dots,n\}$.  We identify the set 
$\{0,1\}^n$  with the set of functions 
$\sigma \colon [n] \to \{0,1\}$
and consider the partial order defined by
$$\sigma \prec \sigma' \iff \sigma(i) \leq \sigma'(i)\ \hbox {\ for all\ } i = 1,
\dots n.$$
Whenever convenient we will identify an element $\sigma \in \{0,1\}^n$ with the subset of $[n]$:
$$ S_\sigma = \{i \in [n] : \sigma(i) = 1\}.$$
Note that a chain 
$\csig = \{\sigma_1  \prec \cdots \prec \sigma_s\}$
in $\{0,1\}^n$ is equivalent to a chain
$\cset = \{S_1 \subset \cdots \subset S_s\}$
of subsets of $[n]$, strictly ordered by inclusion and where  we simply write $S_j$ for $S_{\sigma_j}$.

\begin{definition}\label{def:chain}
A chain $\cset = \{S_1   \subset \cdots \subset S_s\}$ of subsets of $[n]$ is called a chain of $A$-supports if and only if for each $j=1,\dots,s$ 
there exists a vector $v^j \in \ker_\k (A)$ such that
$S_j = \{i \in [n] : v^j_i\not= 0\}$.
\end{definition}
A chain $\cset$ may be alternatively described by the collection 
of mutually disjoint differences: 
$\cdiff = \{D_1,\dots,D_s\}$,
where $D_1 = S_1$ and, for $j \geq 2$, $D_j = S_j \setminus S_{j-1}$.  
Again, we think of $D_j$ as given by the $0,1$ vector with coordinates $1$ in those indices belonging to $D_j$.
When the chain $\cset$ comes from a maximal chain of $A$-supports, we have that $s=m(A)$.  Moreover, as the fact that $A$ is not a pyramid is 
equivalent to the fact that ${\rm ker}_\k(A)$ is not contained in a coordinate hyperplane, the family 
$\cdiff$ is a partition of $[n(A)]$.  A partition of $[n(A)]$ arises from a maximal chain of $A$-supports if its associated chain $\cset$ does.

Recall that given a homogeneous configuration $A$ and a maximal chain $\csig$ of $A$-supports, we have defined the matrix $M_{\csig}(A)$ by:
$$M_{\csig}(A) = \left(\begin{array}{c}A \\\Sigma\end{array}\right)\ ;\quad 
\Sigma = \left(\begin{array}{c}\sigma_1 \\\vdots \\\sigma_{m}\end{array}\right).$$
Note that the matrix $\Sigma$ is row-equivalent to the matrix defined by the partition of differences:
$$ \left(\begin{array}{c}D_1 \\\vdots \\D_{m}\end{array}\right),$$
where the row $D_i$ indicates a row with ones in the positions indicated by $D_i$ and zeroes elsewhere, because the $i$-th row of
this matrix is equal to the difference $\sigma_i - \sigma_{i-1}$ for any $i >1$.

\bigskip

Given a homogeneous non-pyramidal configuration $A$ over $\k$, 
we will denote by $\mathcal I$ the collection of all iterated circuits $I\subset A$ with $d(I)=d(A)$ and 
by $\mathcal M$ the collection of maximal  chains $\csig$ of $A$-supports.  We then have:

\begin{theorem}\label{prop:bij}
For any  homogeneous non-pyramidal configuration $A$ over $\k$, there are maps 
$$
\varphi: \mathcal I \longrightarrow \mathcal M, \quad \psi: \mathcal M \longrightarrow \mathcal I,
$$
such that:
\begin{enumerate}
\item[i)] For  $I \in  \mathcal I$, 
\begin{equation}\label{eq:varphi}
d(I) - \eta(I) = n(A)-1 -{\rm rank}(M_{\varphi(I)}(A)),
\end{equation}
\item[ii)] For  $\csig \in \mathcal M$,
\begin{equation}\label{eq:psi}
d(\psi(\csig)) - \eta(\psi(\csig)) = n(A)-1 -{\rm rank}(M_{\csig}(A)).
\end{equation}
\end{enumerate}
\end{theorem}

\bp
For simplicity of notation we write $d =d(A)$, $n=n(A)$ and $m = m(A)$.
Suppose $I = I_1 \cup \cdots \cup I_p$ is an iterated circuit with $d(I) = d$.  We may assume without loss of generality that 
$ I = \{a_{1},\dots, a_{r }\}.$

We now define $\varphi(I)$ as the chain of $A$-supports whose differences are defined as follows.  Consider the partition
$$D_j = \{i\in [n] : a_i \in I_j\}, \,  j=1,\dots,p, \quad
D_{p+k} = \{r+k\}, \,  k=1,\dots,n  - r.$$
We need to show that this partition comes as the differences of a chain of $A$-supports as in Definition~\ref{def:chain}
and to compute the rank of the corresponding matrix $M_{\varphi(I)}(A)$.

Note that after reordering, we may assume that the matrix defined by the configuration $I$ is affinely equivalent to a matrix of the form
\begin{equation}\label{eq:ic}
 \left(\begin{array}{cccc}
 I_1 & * & \cdots & *\\
0 & \tilde I_2 & \cdots & *\\
 \vdots & \vdots & \ddots & \vdots\\
0&0 & \cdots &\tilde  I_p \\
\end{array}\right),
\end{equation}
where for $j\geq 2$, $I_j$ corresponds to the columns ``containing" $
\tilde I_j$.  

Since $A$ is homogeneous, $I_1$ is a homogeneous circuit and, consequently, the set $D_1$ 
is the support of an element in $\ker_\k(A)$, and this support is minimal. Similarly, 
since the elements in $\tilde I_j$ are a minimally linearly dependent set in $L(A)/L_{j-1}$, it follows that 
$D_1 \cup \cdots \cup D_j$ is the support of an element in  $\ker_\k (A)$.

For $k = 1, \dots, n  - r$, we have that $D_{p+k} = \{r+k\}$.
Note that the assumption $d(I) = d(A) $ implies that 
 that  $a_{r+k}$
is a linear combination of  elements in $I_1 \cup \dots \cup I_p$ for any $k$. As the union of supports is a support,
 we have that $D_1 \cup \dots \cup D_p \cup D_{p+1}\dots \cup D_{p+k}$ 
is the support of some element in $\ker_\k (A)$ for any $k=1, \dots, n-r$.  Then, $\varphi(I)  \in \mathcal M$.  

The matrix $M_{\varphi(I)}(A)$ is row equivalent to a block-upper triangular matrix of the form
\begin{equation}\label{eq:ic2}
 \left(\begin{array}{ccccc}
 I_1 & * & \cdots & * & * \\
  {\mathbf 1} & 0 & \cdots & 0 & 0\\
0 & \tilde I_2 & \cdots & *& * \\
0 & {\mathbf 1} & \cdots & 0 & 0 \\
 \vdots & \vdots & \ddots & \vdots& \vdots\\
0&0 & \cdots &\tilde  I_p& *  \\
0&0 & \cdots &  {\mathbf 1}& 0  \\
0&0 & \cdots &  0& I_{ n-r}  \\\end{array}\right),
\end{equation}
where the symbol ${\mathbf 1}$ represents a row of $1$'s of the appropriate size.

Clearly ${\rm rank}(M_{\csig}(A))$ equals 
$(n -1)$ minus the number of homogeneous configurations $\tilde I_j$, for $j=2,\dots,p$,  
which in turn equals $d(I) - \eta(I)$.  This completes the proof of (\ref{eq:varphi}).

We will now define the map $\psi: \mathcal M \longrightarrow \mathcal I$ and prove (\ref{eq:psi}).
Let $\csig$ be a maximal proper chain of $A$-supports and let us denote by 
$$\cdiff = \{D_1,\dots,D_{m}\}$$
the associated difference sets.   We may assume --up to reordering of $A$-- that the subsets $D_j$ 
are coherently ordered in the sense that for any $k \leq m $, the sets $\{D_1,\dots,D_k\}$ are a partition of $\{1, \dots, |D_1\cup\cdots\cup D_k|\}$.  

Denote by $D_{i_1},\dots,D_{i_p}$ the difference sets with cardinality greater than $1$.  
We must have $D_{i_1} = D_1$ since $A$ is homogeneous.  Let $I_j \subset A$ be the subsets indexed by $D_{i_j}$.  We claim that 
$I = I_1 \cup \cdots \cup I_p$
is an iterated circuit.  

It is clear that $D_1$ equals the non-zero support of $\sigma_1$ and, hence, is a minimally dependent set of $A$,  i.e. a (homogenous) circuit.  
We note next that if $D_s = \{\ell\}$ is a singleton then 
$a_\ell \in L(\bigcup_{1\leq j < s} D_j).$

Let 
$\pi_j \colon L(A) \to L(A)/L_j$ be the projection
where, as before, $L_j = L(I_1 \cup \cdots \cup I_j)$.
Since $D_{j+1}$ is a minimal $A$-support modulo $L_j$, we have that the elements in $\pi_j(I_{j+1})$ 
are minimally linearly dependent.  Moreover, because of the maximality of the chain $\csig$,
if the affine dimension of $\pi_j(I_{j+1})$ is positive, the map $\pi_j$ must be injective on the subset $I_{j+1}$.  Hence, $I$ is an iterated circuit.

Finally, we note that the proof of (\ref{eq:psi}) is completely analogous to that of (\ref{eq:varphi}) 
since after row operation and column and row permutations the matrix $M_{\csig}(A)$ 
will be equal to the matrix (\ref{eq:ic2}) so that the previous computation of its rank applies.
\ep

\begin{remark}
Note that the composition  $\psi\circ \varphi$ of the maps in the statement of Theorem~\ref{prop:bij} is the identity, while
the composition $\varphi \circ \psi$ is the identity up to possible reordering the singletons in the differences of the supports. 
\end{remark}

\smallskip

We are now ready to prove equality~(\ref{eq:itercirc}).

\begin{theorem}\label{th:itercirc}
Let $A \subset \k^{d +1}$ be a homogeneous non-pyramidal
configuration with cardinality $n(A)$ and affine dimension $d(A)$.  Let $\iota(A)$ and  $\rho(A)$ be as in (\ref{eq:iota}) 
and (\ref{eq:rho}), respectively.  Then, as asserted by (\ref{eq:itercirc}),
$$n(A) - 1 - \rho(A) = d(A)  - \iota(A).$$
\end{theorem}

\begin{proof}
Let $I \subset A$ be an iterated circuit 
with $ \iota(A)=\eta(I)$. We may assume without loss of generality that 
$d(I) = d$. Indeed,  by Lemma~\ref{lem:imax} there exists an extension $I'$ of $I$ with $d(I') = d$ but, $\eta(I') \geq \eta(I)$ and, consequently, 
$ \iota(A)=\eta(I')$ as well.
  
Then,  we deduce from~(\ref{eq:varphi}) in Theorem~\ref{prop:bij} that 
$$ n-1 - d  +  \iota(A) ={\rm rank}(M_{\varphi(I)}(A)) \le \rho(A).$$
Conversely, let 
$\csig=\{\sigma_1 \prec \dots \prec \sigma_{m }\}$ 
be  a maximal  chain of $A$-supports and assume that 
${\rm rank}(M_{\csig}(A)) = \rho(A)$. Then,  using~(\ref{eq:psi}) in Theorem~\ref{prop:bij} we have:
$$\rho(A) = n-1 - (d(\psi(\csig)) - \eta(\psi(\csig)) )
= n-1 - d + \eta(I) \le n-1 - d  +  \iota(A),$$
which ends the proof.
\end{proof}

\section{Gale Dual: Proof of (\ref{eq:flags})}\label{sec:right=}

In this section we will recall the construction of the $\lambda$-invariant introduced in \cite{CC} and prove  identity (\ref{eq:flags}) in Theorem~\ref{th:gale}. 
This gives a formula for the dual defect of a complex toric variety $X_A$ in terms of a Gale dual of $A$.  
In the final subsection, we describe the behavior of  all the invariants in case $A$ is a pyramid. 
As in the previous section, we will work with configurations defined over an arbitrary field $\k$ of characteristic zero.

Let $A = \{a_1,\dots,a_n\} \subset \k^{d + 1}$ be a configuration such that $\kappa(A) = n -d -1$ is the dimension of $W=\ker_\k(A)$.  
Note that if $d(A)=d$, in case $A$ is homogeneous $\kappa(A)$ coincides with $m(A)$, which is the dimension of the space of affine relations of $A$, and $\kappa(A) = m(A)-1$ otherwise.

We have a short exact sequence
\beq\label{eq:exseq}
\begin{tikzcd}
0 \arrow[r] & W \arrow[r,"\iota"] & \k^n \arrow[r,"\alpha"] & \k^{d+1} \arrow[r] & 0,
\end{tikzcd}
\eeq
where $\alpha$ maps the standard basis element $e_i \in \k^n$ to 
$a_i \in \Z^{d+1}$.   Dualizing (\ref{eq:exseq}) we get
\beq\label{eq:exseqdual}
\begin{tikzcd}
0 \arrow[r] & (\k^{d+1})^* \arrow[r, "\alpha^*"] & (\k^n)^* \arrow[r,"\iota^*"] & W^* \arrow[r] & 0.
\end{tikzcd}
\eeq
Choosing a basis we may identify $W^* \cong \k^{\kappa(A)}$ and we denote by $b_i := \beta(\xi_i)\in \k^{\kappa(A)}$, where $\xi_i$ is the standard dual basis 
of $(\k^n)^*$.  We denote by $B$ the $n\times \kappa(A)$ matrix whose rows are the $b_i$'s and note that the columns of 
$B$ are a basis of $\ker_\k A$.  We will refer to the configuration 
$\{b_1,\dots,b_n\}$ as a Gale dual of $A$.  

In what follows it will be useful to consider the linear matroid of rank 
 $\kappa(A)$ defined by the configuration
$ \{b_1,\dots,b_n\}$ on $\{1,\dots,n\}$.  Note that $A$ is non-pyramidal if and only if none of the $b_j=0$; i.e. if and only if the matroid $B$ has no loops.

Given any subset $C \subset B$ we denote
\begin{equation}\label{eq:sum} 
\fs(C) = \sum_{b\in C} b.
\end{equation}
 The configuration $A$ is homogeneous if and only if any Gale dual $B$ satisfies the {\em dual homogeneity} condition
$$ \fs(B)= \sum_{b\in B} b \ =\ 0.$$
In this case we refer to $B$ as a {\em dual-homogeneous} configuration.

Given a configuration $A$ we have defined in~(\ref{eq:overline}) and~(\ref{def:Chom}) two associated homogeneous configurations $\bar{A}$ and $A^h$.  
It is easy to see that their corresponding Gale duals may be obtained from a Gale dual $B$ of $A$ as follows:  A Gale dual $\bar{B}$ of $\bar{A}$ is the 
$m(A)$-dimensional configuration obtained by projecting $B$ to the quotient $L(B)/L(\fs(B))$, and a Gale dual of $A^h$ is the dual-homogeneous configuration of rank $m(A) + 1$
\beq\label{eq:homogB}
B^H = \{-\fs(B)\} \cup B .
\eeq

A subset $F\subset B$ is called a flat  if and only if
$ F = L(F) \cap B$.  
The {\em rank} of a flat $F$ is defined by $\rank(F) = \dim L(F)$.
A flag in $B$  is a collection of flats
$\mathcal F = \{F_1 \subset \cdots \subset F_\ell\},$
with $\rank F_j = j$ for every $j=1, \dots, \ell$. 
We call $\ell$ the length of the flag.

\begin{definition}\label{def:nsf}
We say that a flag $\mathcal F=  \{F_1 \subset \cdots \subset F_\ell\}$  in $B$ is
 {\em non-splitting} if 
$$\fs(F_1)\not=0, \hbox{\  and \ } \fs(F_j)  \not\in L(F_{j-1}), \hbox{\  for any\  } 2\leq j \leq \ell .$$
%
Given a vector configuration $B$ defined over $\k$ we denote by 
$\lambda(B)$ the maximal length of a non-splitting flag $\mathcal F$ in $B$.
  We carry the notion to the $A$-side by setting $\lambda(A) = \lambda(B)$ for any Gale dual $B$ of $A$.  
\end{definition}
  
Since two Gale duals of a given 
  homogeneous configuration $A$ differ by the action of an element in $GL(m,\k)$, it is clear that $\lambda(A)$ does not depend on the choice of a Gale dual. 
Note that if $A$ is homogeneous and $B$ is a Gale dual of $A$ then 
$\fs(B)=0$ and therefore $\lambda(A) \leq m(A)-1$.

\smallskip

Before stating our main result in this section Theorem~\ref{th:gale}, we collect some preliminary results.

\begin{lemma}\label{circuit-lemma}
Let $A= \{a_1, \dots, a_{n}\}$ be a homogeneous, non-pyramidal  configuration over $\k$ and suppose that
$Z \subset A$ is a circuit.  Set
$A_1 = (L(Z)  \cap A) \setminus Z$
and $A_2 = A \setminus (Z \cup A_1)$. Let $C = \{b_i, \, a_i \in Z\}$, $B_1= \{b_i, \, a_i \in A_1\}$, $B_2= \{b_i, \, a_i \in A_2\}$ denote the 
corresponding collections in a Gale dual configuration $B =\{b_1, \dots, b_{n}\}$ of $A$.  Then
\begin{enumerate}
\item[i)] The set $B \setminus C = B_1 \cup B_2$ is a codimension-one flat in $B$.
\item[ii)] The elements in $B_1$ are linearly independent.
\item[iii)] The set $B_2$ is a flat of codimension $n(A_1)+1$.
\item[iv)] $ L(B_1 \cup B_2) = L(B_1) \oplus L(B_2).$
\end{enumerate}
\end{lemma}

\begin{proof} Since $Z$ is a circuit, its Gale dual is one-dimensional.  But the Gale dual of $Z$ is the projection of $C$ to the quotient
$L(B)/L(B_1 \cup B_2)$.  Hence, $L(B_1 \cup B_2)$ has codimension one.  
Moreover, since $Z$ is a circuit it is non-pyramidal and hence no element of $Z$ projects to zero on the quotient.  Hence 
$B_1 \cup B_2$ is a flat.

Since $d(Z) = d(Z \cup A_1)$ we must have 
$$\dim(L(B_1 \cup B_2)) = \dim(L(B_2)) + n(A_1).$$
Hence the elements in $B_1$ are linearly independent.

Similarly, $L(B_2)$ is a subspace of codimension
equal to the dimension of the quotient $ L(B) / L(B_2)$. This is the underlying space of a 
Gale dual of $ Z \cup A_1$, and therefore has codimension $n(A_1)+1$.  Moreover since 
$L(Z) = L(Z \cup A_1)$ and $Z$ is not pyramidal, then $A_1$ is not a pyramid either and, consequently, $B_2$ is a flat.
Item iv) follows from items ii) and iii). 
\end{proof}

\begin{lemma}\label{quotient-flag}
Let $B =\{b_1, \dots, b_n\}$ be a dual-homogeneous  configuration over $\k$  
and let $\lambda$ be the maximal length of a non-splitting flag in $B$.  Let 
$F_1 \subset \cdots \subset F_{\lambda} \subset B$
be a non-splitting flag.  
Set 
$\tilde{B}= \Pi(B \setminus F_1)$,
where $\Pi \colon L(B) \to L(B) / L(F_1) $ is the canonical projection. 
Let $\tilde \lambda$ denote the maximal length of a non-splitting flag in $\tilde B$.
Then 
$\tilde \lambda = \lambda -1$.
\end{lemma}

\begin{proof}  
Suppose 
$H_1 \subset \cdots \subset H_{\tilde\lambda}$
is a maximal non-splitting flag in the dual-homo\-ge\-neous configuration $\tilde{B}$.  For $i=2,\dots,\tilde\lambda+1$, let 
$G_i$ denote the flat of $B$ defined by 
$G_i =\{b\in B : \Pi(b) \in H_{i-1}\}$.
It follows that
$F_1 \subset G_2 \subset \cdots \subset G_{\tilde \lambda + 1}$ 
is a non-splitting flag in $B$.  Hence 
$\lambda \geq \tilde\lambda+1$.

On the other hand, it is clear that the projections
$\Pi(F_2) \subset \cdots \subset \Pi(F_\lambda)$
define a non-splitting flag in $\tilde{B}$ and therefore 
$\tilde \lambda \geq \lambda -1$.
\end{proof}

We now recall for later use the statement of Lemma~22  in~\cite{CC}.

\begin{lemma}\label{lemma22}
Let $B$ be a dual-homogeneous configuration and $\Lambda \subset L(B)$ a line.  If 
 $B$ has a non-splitting flag $\mathcal F$ of length $\ell$, then $B$ has a non-splitting flag $\mathcal G$ of length $\ell$ such that 
 $\Lambda \cap L(G_\ell) =\{0\}$.
\end{lemma}

Next we will consider the case of maximal non-splitting flags in non-homogeneous configurations.

\begin{lemma}\label{lambda-nothomog}
Let  $B = \{b_1,\dots, b_n\}$
be a configuration and suppose $\fs(B) \not =0$.  Let $B^H$ denote its homogenization; i.e. the configuration
$B^H = \{-\fs(B),b_1,\dots, b_n\}$.
Then
$$\lambda(B) = \lambda(B^H) + 1,$$
where $\lambda(B)$ and $\lambda(B^H)$ are
the maximal lengths of non-splitting flags in $B$ and $B^H$.
\end{lemma}

\begin{proof}
Set $\lambda=\lambda(B)$ and $\lambda^H = \lambda(B^H)$.
Suppose $\hat F_1 \subset \cdots \subset \hat F_{\lambda^H}$
is a maximal non-splitting flag in $B^H$.  By Lemma~\ref{lemma22}, we may assume that 
$\fs(B) \not \in  L(\hat F_{\lambda^H})$ and therefore all $\hat F_i$ are flats in $B$.  
We simbolize this by dropping the $\ \hat{}\ $ from the notation.

Suppose now that $F_1 \subset \cdots \subset F_{\lambda^H}$ 
may not be extended to a longer non-splitting flag in $B$.  
Let $G_1,\dots,G_r$ denote the flats of rank 
$\lambda^H + 1$ in $B$ containing $F_{\lambda^H}$.  
Then, for each $G_i$ we must have $\fs(G_i) \in  L(F_{\lambda^H})$.  
But this implies
$$\fs(B) = \sum_{i=1}^r \fs(G_i) - (r-1) \fs(F_{\lambda^H}) \in  L(F_{\lambda^H}),$$
which is a contradiction.  Hence, $\lambda(B) \geq \lambda(B^H) + 1$.

To prove the opposite inequality, we begin by observing that given a maximal non-splitting flag $\mathcal F =\{F_1 \subset \cdots \subset F_{\lambda}\}$  in $B$,  
there exists $j=1,\dots,\lambda$ such that $\fs(B) \in L(F_j)$.  Otherwise, we could regard $\mathcal F$ as a non-splitting flag 
of length $\lambda$ in $B^H$ and, consequently, $\lambda(B^H) \geq \lambda(B)$ which would contradict the inequality already proven.  
Thus, let $j = 1,\dots, \lambda$ be such that 
$$\fs(B) \in L(F_j) \setminus L(F_{j-1}).$$
If $j = \lambda$, then $F_1 \subset \cdots \subset F_{\lambda-1}$ may be regarded as a non-splitting flag in $B^H$ implying that $\lambda^H \geq \lambda - 1$ as desired.

Thus, it suffices to show that we can always find a non-splitting flag $\mathcal G$ of length $\lambda$ in $B$ so that 
\begin{equation}\label{top}\fs(B) \in L(G_\lambda) \setminus L(G_{\lambda-1}).
\end{equation}
  Suppose then that $\fs(B) \in L(F_j) \setminus L(F_{j-1})$ and $j < \lambda$.  Let $F_j,H_1,\dots,H_r$ denote the distinct $B$-flats of rank $j$ in $F_{j+1}$
   containing $F_{j-1}$.  There must exist an index $i=1,\dots,r$, such that 
\begin{equation}\label{spl}
\fs(H_i) \not \in L(F_{j-1})
\end{equation}
since, otherwise,
$$\fs(F_{j+1}) = \fs(F_j) + \sum_{i=1}^r \fs(H_i) - r\,\fs(F_{j-1}) \in L(F_j)$$
which is not possible since $\mathcal F$ is non-splitting. 

We next claim that the index  $i=1,\dots,r$ may be chosen so that  (\ref{spl}) is satisfied  and  $\fs(F_{j+1}) \not \in L(H_i)$.  
Indeed, suppose there exists a unique flat $H_i$ satisfying (\ref{spl}).  Then
$$\fs(F_{j+1}) \equiv \fs(F_j) + \fs(H_i)\ \hbox{mod} 
L(F_{j-1})$$
but then $\fs(F_{j+1}) \not \in L(H_i)$ since $\fs(F_j)$ and $ \fs(H_i)$ are linearly independent modulo $L(F_{j-1})$.  On the other hand if there is more than 
one $H_i$ satisfying (\ref{spl}) we can certainly choose one of them whose span does not contain $\fs(L_{j+1})$.  

Thus, we can construct a new non-splitting flag $\mathcal G$ replacing $F_j$ by $H_i$ so that $\fs(B) \in L(G_{j+1}) - L(G_j)$.  Continuing inductively we construct a flag satifying (\ref{top}).
\end{proof}

The following simple example illustrates these invariants in the non-homogeneous case.

\begin{example}\label{ex:123}
Let $A$ be the non-homogeneous configuration 
$$A = [1\ 2\ 3] \subset \Q^{1\times 3}.$$  Then a Gale dual $B \subset \Q^2$ and its homogenization $B^H$ are  given by the rows of the matrices:
$$B = \left(\begin{array}{rr}-2 & 0 \\1 & -3 \\0 & 2\end{array}\right)\ ;\quad
B^H = \left(\begin{array}{rr}1 & 1 \\-2 & 0 \\1 & -3 \\0 & 2\end{array}\right)\,.$$
Clearly, $\lambda(B) = 2$ and $\lambda(B^H) = 1$.  Note that a Gale dual of 
 $\bar{A}$ is the one-dimensional configuration 
$\bar{B} = \{1,-2,1\} \subset \Q$ and $\lambda(\bar{B}) = 0$.
\end{example}

Finally, we need a result, somewhat analogous to the previous Lemma, 
about the behavior of the $\rho$-invariant under homogenization.  

\begin{lemma}\label{rho-nothomog}
Let $A$ be a non-pyramidal, non-homogeneous configuration with cardinality $n$ and let $A^h$ be as in (\ref{def:Chom}).  Suppose
$\csig = \{\sigma_1 \prec \cdots \prec \sigma_{m+1}\}$
is a maximal chain of supports of elements in $\ker_\k(A)$.  Let $\Sigma$ be the 
matrix whose rows are the elements $\sigma_1,\dots,\sigma_{m+1}$.  Then
$$\rho(A^h) \geq \rank \left(\begin{array}{c}A \\\Sigma\end{array}\right).$$
\end{lemma}

\bp
For each $j = 1,\dots,  {m+1}$, let 
$v_j  \in \ker_\k (A)$ be such that ${\rm supp}(v_j) = \sigma_j$.  Then the vector
$v^h_j = (-\sum_i\,v^i_j, v_j) \in \ker_\k (A^h)$, where $v^i_j$ are the components of $v_j$. As $A$ is non-homogeneous, 
there is a minimal index $j$ for which $\sigma(v_j)\not=0$. After adding a
multiple of $v_j$ to  any $v_k$ with $k > j$, if necessary, we can  assume that  
the supports $\sigma_j^h={\rm supp}(v^h_j)$, $j=1, \dots, m+1=m(A^h)$
define a maximal chain ${\csig}^h$ of $A^h$-supports. 

Note  that $A^h$ is not a pyramid 
because $A$ is not a pyramid and it is non-homogeneous.  Denote by $\Sigma^h$ the matrix with rows the elements of $\csig^h$.  Then,
$$\left(\begin{array}{c}A^h \\\Sigma^h\end{array}\right) = \left(\begin{array}{cc}1 & \mathbf{1_{n}}\\0 & A \\ \ast & \Sigma\end{array}\right).$$
Hence
$$\rho(A^h) \geq \rank \left(\begin{array}{c}A^h \\\Sigma^h\end{array}\right) \geq \rank
 \left(\begin{array}{c}\mathbf{1_n}\\ A \\ \Sigma\end{array}\right) = \rank \left(\begin{array}{c} A \\ \Sigma\end{array}\right),$$
where the last equality follows from the fact that
$\sigma_{m(A^h)}^h = \mathbf{1_{n+1}}$ is
the all $1$-vector since $A^h$ is homogeneous.
\ep

We are now ready to prove equality~\ref{eq:flags}.

\begin{theorem}\label{th:gale}
Let $A \subset \k^{d +1}$ be a homogeneous, non-pyramidal, configuration with cardinality $n$ and affine dimension $d$.  
Let $m = n - d -1$ and let $\rho(A)$ be as in (\ref{eq:rho}).  Then, as asserted by (\ref{eq:flags})
$$n - 1 - \rho(A) = m -1 - \lambda(A).$$
\end{theorem}

\begin{proof}
We proceed by induction on $m$. If $m=1$ 
then $\lambda(A) = 0$ and $m - 1 - \lambda(A) =0$.  On the other hand, since $A$ is a circuit, there are no proper supports 
in ${\mathcal S}(A)$ and therefore $\rho(A) = n-1$.  Thus, (\ref{eq:flags}) holds.

We assume that the equality holds for configurations of codimension less than $m$ and begin by showing that 
\begin{equation}\label{first_ineq} 
n-\rho(A) \geq m - \lambda(A).
\end{equation}

  Let $\csig = \{\sigma_1 \prec \cdots \prec \sigma_{m}\}$
be a maximal  chain of $A$-supports such that $\rho(A) = \rank(M_{\csig}(A))$.  
We may assume without loss of generality that 
$$\sigma_1 = (\overbrace{1,\dots,1}^{s},0,\dots,0).$$
Then, since $\sigma_1$ is a minimal support in $\ker_\k(A)$,  $Z=\{a_1,\dots,a_s\}$ is a (homogeneous) circuit  in $A$.
In particular, $d(Z) = s -2 $ and, after affine transformation, we may assume that $  L(Z) $ 
is the subspace of $\k^{d+1}$ spanned by $e_1,\dots,e_{s-1}$.  

Let $A_1= (A \cap L (Z))  \setminus Z$ and assume that (if non-empty) $A_1= \{a_{s+1},\dots,a_{s+k}\}$.
Thus, we may write 
$$A = \left(\begin{array}{ccc}Z & A_1 & \ast \\0 & 0& A_{2}\end{array}\right),$$
where no column of $A_2$ is zero.
By Lemma~\ref{circuit-lemma} we have that the vectors
 $\{b_{s+1},\dots, b_n\}$ define a codimension-one flat in $B$, while the vectors $\{b_{s+1},\dots,b_{s+k}\}$ are linearly independent. Moreover,
$ L(\{b_{s+1},\dots,b_n\}) = L(\{b_{s+1},\dots,b_{s+k}\}) \oplus L(\{b_{s+k+1},\dots,b_n\})
$.

 Thus,  $B$ can be assumed to be of the form
$$B = \left(\begin{array}{ccc}v & \ast & \ast \\0 & E_k & 0 \\0 & 0 & B_2 \end{array}\right),$$
where $v$ is a column vector of length $s$ with support $\sigma_1$ (a generator of $\ker_\k(Z)$), 
$E_k$ is a $k\times k$ identity matrix, and the columns of $B_2$ are a basis of $\ker_\k(A_2)$.  Note that
since $A$ is non-pyramidal, no row of $B_2$ can be zero and that  the last $n-s-k$ coordinates
of any vector in $\ker_\k(A)$ is a vector in $\ker_\k(A_2)$. 

Now, for any $j \in \{1, \dots, k\}$, the element
$a_{s+j}$ is in the  span of $Z$, and therefore there exists a vector $v_j \in \ker_\k(A)$ with
support $\{1, \dots, s\} \cup \{s+j\}$.  Hence  any {\em maximal} chain $\csig$ should satisfy that
 there exists an index $\ell(j)$ such that  the difference of the supports
 ${\rm supp}(\sigma_{\ell(j)}) \setminus {\rm supp}(\sigma_{\ell(j)-1})$ equals
 $\{s+j\}$. Otherwise, we would be able to add another support to the chain, contradicting the maximality.

Then, $M_{\csig}(A)$ is then row equivalent to a matrix of the form
$$\left(\begin{array}{ccc}Z & A_1 & \ast \\{\mathbf 1_s} & 0 & 0 \\0 & E_k & 0 \\0 & 0 & A_{2} \\0 & 0 & \tilde\Sigma\end{array}\right),$$
where as before ${\mathbf 1_s}$ is a row vector of  length $s$ whose entries are all $1$, and 
$\tilde \Sigma$ is a matrix whose rows define a maximal  chain of  $A_2$-supports.
Therefore, since $Z$ is homogeneous we have that
\begin{equation}\label{eqn:rho}
\rho(A) = \rank(M_{\csig}(A)) = (s-1) + k + \rank \left(\begin{array}{c}A_{2} \\ \tilde\Sigma \end{array}\right).
\end{equation}

We now consider two cases. Assume first  that $A_{2}$ is  
homogeneous and so the rows of $B_2$ form a dual-homogeneous configuration
which is a Gale dual of $A_{2}$. Then,
$$\rank \left(\begin{array}{c}A_{2} \\ \tilde\Sigma \end{array}\right) \leq \rho(A_2),$$
and consequently
\[\rho(A) \leq (s-1) + k + \rho(A_{2}).\]
We have $n(A_2) = n -s -k$, $m(A_{2}) = m-1-k$, and, by inductive hypothesis,
$$n(A_2) - \rho(A_2) = m(A_2) - \lambda(A_2).$$
Hence
\beas
n - \rho(A) & \geq & n - s  - k + 1 - \rho(A_2)
= n(A_2) - \rho(A_2) + 1\\
&=& m(A_2) - \lambda(A_2) + 1
= (m-1-k) - \lambda(A_2) + 1\\
&=& m - (\lambda(A_2) + k)
\geq  m- \lambda(A),
\eeas
where the inequality $k + \lambda (A_2) \leq \lambda (A)$ follows from the fact
that every non-splitting flag of length $\ell$ in $B_2$ may be extended to a non-splitting flag of length $\ell + k$ 
by adjoining the $k$ linearly independent elements of $B_1$ one at a time.

\medskip

Suppose now that $\fs(B_2) \not= 0$.  Let $ B_2^H$ be the homogenization of $B_2$.  Then,  $B_2^H$ is a Gale dual of
$A_{2}^h$ and 
$$m(A_{2}^h)= m(A_2), \quad n(A_{2}^h)= n(A_2)+1= n-s-k+1.$$
It follows from Lemma~\ref{rho-nothomog} that 
$$\rho(A_2^h) \geq \rank \left(\begin{array}{c}A_{2} \\ \tilde\Sigma \end{array}\right)$$ 
and therefore we may argue as in the previous case and have:
\beas
n - \rho(A) & \geq & n - s  - k + 1 - \rho(A^h_2)
= n(A^h_2) - \rho(A^h_2) \\
&=& m(A^h_2) - \lambda(A^h_2) 
= (m-1-k) - \lambda(A_2) \\
&=& m - (\lambda(A^h_2) + k + 1)
\geq  m- (\lambda(A_2)+k)\\
&\geq & m-\lambda(A),
\eeas
where the next to last inequality follows from Lemma~\ref{lambda-nothomog} and the last inequality follows, as above,  from the fact
that every non-splitting flag of length $\ell$ in $B_2$ may be extended to a non-splitting flag of
length $\ell + k$ by adjoining the $k$ linearly independent elements of $B_1$ one at a time.
 This completes the proof of (\ref{first_ineq}).

\bigskip

To complete the proof of Theorem~\ref{th:gale} we now show that 
\begin{equation}\label{second_ineq} n-\rho(A) \leq m - \lambda(A).
\end{equation}
  Set $\lambda = \lambda(A)$ and suppose that 
$\ F_1 \subset \cdots \subset F_\lambda\ $
is a maximal non-splitting flag in $B$.  Let $\tilde B$ denote the projection of $B \setminus F_1$ to the quotient 
$L(B)/ L(F_1)$.  Then $\tilde B$ may be seen as the Gale dual of the (homogeneous) configuration $A_1 = \{a_i \in A : b_i \not \in F_1\}$.
We assume without loss of generality that $|F_1| = s$,
 and that $F_1$ consists of the last $s$ vectors in $B$. So,  
$A_1$ consists of the first $n-s$ vectors in $A$. As $m(A_1)=m-1$ and $n(A_1)=n-s$,  we have that $d(A_1)+1=n-s-m+1$. 
We can assume that $A$ is of the form
$$A = \left(\begin{array}{cc}A'_1 & \ast \\0 & Z\end{array}\right),$$
where $Z \in \k^{(s-1) \times s}$.

From Lemma~\ref{quotient-flag} we have that $\lambda(\tilde B)  = \lambda(B) -1$.  Thus, by inductive assumption we have that
$n - s- \rho(A'_1) = (m - 1) - \lambda(\tilde B_1) = m - \lambda(A)$,
and therefore 
\begin{equation}\label{eq:rhoA1}
\rho(A'_1) + s = n - m + \lambda(A).\end{equation}
Suppose now that $\tilde\sigma_1 \prec \cdots \prec \tilde\sigma_{m-1}$ is a maximal chain of $\ker_\k(A'_1)$-supports 
computing $\rho(A'_1)$.  Set
$$\sigma_j = (\tilde\sigma_j,\overbrace{0,\dots,0}^r)\ ; \quad j=1,\dots,m-1,$$
and $\sigma_{m} = (1, \dots, 1) = {\mathbf 1}_n$.
Then, $\sigma_1 \prec \cdots \prec \sigma_{m}$ is a maximal chain of supports for $\ker_\k(A)$ and we have
$\rho(A) \geq \rank(M_{\csig}(A))$.
But $M_{\csig}(A)$ is row-equivalent to the matrix
$$\left(\begin{array}{cc}A'_1& \ast \\\tilde\Sigma & 0 \\0 & Z \\0 & {\mathbf 1}_s\end{array}\right)$$
whose rank is $\rho(A'_1) + s$ since $F_1$ is a non-splitting flat which implies that $Z$ is not homogeneous; i.e. 
${\mathbf 1}_s$ is not in the rowspan of $Z$.  Hence, we deduce
from~(\ref{eq:rhoA1}) the inequality
$$
\rho(A) \geq \rho(A'_1) + s = n - m + \lambda(A),
$$
as wanted. 
\end{proof}

\bigskip

\subsection{\bf {Pyramids}} \label{sec:pyramid}

 Let us consider the invariants defined in the previous sections in the pyramidal case.

  \begin{definition}\label{def:pyr}
 Let $A$ be a  configuration over a field $\k$ of characteristic zero.  We say that 
 $A$ is a pyramid if all points in $A$ but one lie on an (affine) hyperplane.
When $A$ is not a pyramid,  we define $p(A)=0$. If $A$ is a pyramid, let $a \in A$
such that the points in $A \setminus\{a\}$ lie in a hyperplane not contaning
$a$. Then, $p(A)$ is defined, inductively, as $p(A \setminus\{a\})+1$.  We will call $p(A)$ the pyramidal index of $A$.
\end{definition}

Alternatively, we may describe the quantity $p(\bar{A})= p(A)$ as the 
 number of indices $i$ with $b_i=0$ in a Gale dual  configuration $B$ of $\bar{A}$.

Assume $A$ is a pyramid, that is $p(A) \ge 1$. The matrix $A$  has the following shape, after
 left multiplication by an invertible matrix and probably renumbering the elements of $A$:
 \begin{equation}
\label{eq:pyr}
\left(\begin{array}{cc}
 E_{p(A)} & 0 \\
0 &  A'
\end{array}\right),
\end{equation}
 where $E_{p(A)}$ is the $p(A)\times p(A)$ identity matrix and  $p(A')=0$.
Note that  $A$ is homogeneous if and only if $A'$ is homogeneous.
 
 The following summarizes the behavior of the invariants studied above in the pyramidal case.  
 The proofs are straightforward and are left to the reader:

\begin{theorem}\label{th:pyramid} 
Let $A$ be a configuration over $\k$ with pyramidal index $p(A)$.  Let $A'$ be as in (\ref{eq:pyr}).  Then
\begin{enumerate}
\item[i)] $n(A) = n(A') + p(A)$;
\item[ii)] $d(A) = d(A') + p(A)$;
\item[iii)] $m(A) = m(A')$;
\item[iv)] $\iota(A) = \iota(A')$;
\item[v)] $\rho(A) = \rho(A') + p(A)$;
\item[vi)] $\lambda(A) = \lambda(A')$
\end{enumerate}
\end{theorem}

It is then clear that both (\ref{eq:itercirc}) and (\ref{eq:flags})  fail if $p(A)>0$.  
On the other hand, in the geometric case, 
we have that
$${\mathrm {def}}(X_A) = {\mathrm {def}}(X'_A) + p(A).$$
Thus we see that only the quantity
$$d(A) - \iota(A)$$
computes the geometric defect; while the $m(A)-1-\lambda(A)$ and $n(A)-1-\rho(A)$ invariants only give the right answer if $A$ is non-pyramidal.


\section{Defect and Cayley Decompositions}\label{sec:cayley}

Furukawa and Ito  introduced in~\cite{FI}
an alternative method for computing the dual defect of a projective toric variety $X_A$ over $\C$
in terms of Cayley decompositions of the configuration $A$. In this section, we translate their
results to the Gale dual setting and describe some applications.
We begin with a brief description of  Furukawa and Ito's approach, referring to~\cite{FI} for details and proofs.   Our translation
to the Gale dual setting starts with Proposition~\ref{lem:ABC} in \S~\ref{ssec:Gale} and yields
Theorem~\ref{th:fi_dual}, which is a vast improvement over \cite[Theorem~25]{CC}. 
This makes it possible to give a detailed description of dual defective toric varieties.
We conclude this section by sketching some of its consequences in \S~\ref{ssec:BFI}.

\medskip

Since in this section we are interested in the geometric setting, we will assume from the start that 
$A \subset  \Z^{d+1}$ is a homogeneous configuration, which arises as $A = \overline{A'}$ with
$A'\subset \Z^d$ as in~(\ref{eq:overline}).
Following ~\cite{FI}, we denote by $\langle A' - A'\rangle$ the affine span of $A'$ over $\Z$ and 
assume that $\langle A' - A'\rangle=\Z^d$ and, consequently, $d(A) = d$.
As before, let $n = n(A)$ denote the cardinality of $A$ and $m = m(A) = n - d - 1$ the rank of a Gale dual $B$ of $A$.  

\begin{definition}\label{def:cayley}
We say that a homogeneous integer configuration $A$ has a {\em Cayley decomposition} of length $r$ if $A$ is  equivalent, up to a $\Z$-linear isomorphism, to a configuration of the form:  
\begin{equation}\label{eq:Cayley}
 \{e_0\} \times A_0 \cup \{e_1\} \times A_1 \cup \dots \cup
\{e_r\} \times A_r,
\end{equation}
where 
$A_0,\dots,A_r \subset \Z^{d-r}$ are non-empty finite configurations and $e_0, e_1, \dots, e_r$ is the canonical basis in $\Z^{r+1}$. We write
\begin{equation}\label{eq:Cayley_notation}
A = A_0 \ast \cdots \ast A_r.
\eeq
\end{definition}
Note that we may view any homogeneous configuration $A$ as admitting a Cayley decomposition of length $0$.

If $A = A_0 \ast \cdots \ast A_r$ then, up to left multiplying by an invertible integer matrix and permuting the columns, we can write the associated matrix $A$ in the form:
\begin{equation}\label{eq:Cayley_matrix}
\left(\begin{array}{cccc}
{\mathbf 1} & {\mathbf 0} & \cdots & {\mathbf 0} \\
{\mathbf 0} & {\mathbf 1} & \cdots & {\mathbf 0} \\
0 & 0 & \ddots &0 \\
{\mathbf 0} & {\mathbf 0} & \cdots & {\mathbf 1} \\
A_0 & A_1 & \cdots & A_r
\end{array}\right),
\end{equation}
where ${\mathbf 1}$ (resp. ${\mathbf 0}$)  is used to denote row vectors (of appropiate lengths) with all entries equal to $1$ (resp. $0$).

The first link between dual defect and Cayley decompositions appears in the work of Di~Rocco~\cite{Di}, where it is shown that if $X_A$ is a defective 
{\em smooth} toric variety, then $A$ admits a Cayley decomposition of length $r = (d + n)/2$ and such that the convex hull of the $A_i$'s have the same normal fan.  
Similar results were obtained by Casagrande and Di Rocco~\cite{CDi} for normal $\Q$-factorial varieties.  More generally, it is shown in~\cite{CC,E1}
 that every dual defective toric variety admits a Cayley decomposition of positive length.  This result was further extended by Ito~\cite{I}, who proved that every dual defective
  toric variety admits a Cayley decomposition of length at least $\defect(X_A)$.

Before we can state the main result of \cite{FI} computing the dual defect in terms of Cayley decompositions, 
we need to introduce a further invariant of a Cayley decomposition.

\begin{definition}\label{def:join}
We say that a Cayley decomposition as above is a {\em join}  
if the  map $(m_0, \dots, m_r) \mapsto m_0 + \dots m_r$  defines an isomorphism
$$\langle A_0 - A_0\rangle \oplus \cdots \oplus \langle A_r - A_r\rangle \to  \Z^{d-r}.$$
\end{definition}

The following is a restatement of Theorem~1.3 and Corollary~1.5 in \cite{FI}:

\begin{theorem}[Furukawa-Ito]\label{th:fi}.  Let $A = \overline{A'} \subset \Z^{d+1}$ be a homogeneous configuration such that $\langle A' - A' \rangle = \Z^d$.  Then, there exist
\begin{itemize}
\item An integer $r, \, 0 \le r \le d$,  and a Cayley decomposition $A = A_0 \ast \cdots \ast A_r$, 
\item An integer $c, \,  0 \le c\le d-r$, and a subgroup $S \subset \Z^{d+1}$ with $\rank(S) =c$,
\end{itemize}
such that if  $\pi \colon \Z^{d+1} \to \Z^{d+1}/S \simeq \Z^{d+1-c}$ is the natural projection, then
\begin{enumerate}
\item[i)] $\pi(A) = \pi(A_0) \ast \cdots \ast \pi(A_r)$ is a join (with $d(\pi(A)) = d-c$) and
\item[ii)] $\defect(X_A) = r - c$.
\end{enumerate}
This Cayley decomposition is unique in the following sense: given any other Cayley decomposition $A = A'_0 \ast \cdots \ast A'_{r'}$ and subgroup $S'$ satisfying i) and ii) above, we have $S \subset S'$ and for each $j=0,\dots,r$ we have
$$A_j = A'_{i_1} \ast \cdots \ast A'_{i_s}$$
for some $A'_{i_k}$. In particular, $r \le r'$ and $c \le c'$.
\end{theorem}
 
 It is convenient to introduce the affine invariant $\theta(A)$.

\begin{definition}\label{def:theta}
Given a Cayley decomposition $A = A_0\ast A_1 \ast \cdots \ast A_r$  and  a subgroup $S$ of rank $c$
 verifying  item i) in Theorem~\ref{th:fi},   we call $\theta(A)$ the maximum value of the differences $r- c$ as we run over all Cayley decompositions and subgroups.

We say that  a given Cayley decomposition and subgroup is a $\theta$-decomposition of $A$ if $r - c = \theta(A)$.  
Furthermore, we say that it is an {\em FI-decomposition} of $A$ if it is a $\theta$-decomposition and its length $r$ is the smallest possible for a $\theta$-decomposition.
\end{definition}

Note that there might be different Cayley decompositions of a configuration $A$ which are incomparable (see e.g. the prism configuration in Example~\ref{ex:return}.) Theorem~\ref{th:fi} asserts that  given an integer configuration $A$ as above, there exist a unique FI-decomposition with $r-c= \theta(A)$ with $r$ minimal among all $\theta$-decompositions of $A$.  

Moreover, in the geometric case, we have that $\defect(X_A) = \theta(A)$.
On the other hand, it follows from (\ref{eq:defect_rho}) that if $A$  is a homogeneous, non-pyramidal, configuration 
$A= \overline{A'}\subset  \Z^{d+1}$ satisfying $\langle A' - A' \rangle = \Z^d$, then
$\defect(X_A) = n(A) - 1 - \rho(A)$.  Hence, together with 
Theorems~\ref{th:itercirc} and~\ref{th:gale} we have that for such a configuration $A$:
\beq\label{eq:equalities}
\theta(A) = m(A) - 1 - \lambda(A) = n(A) - 1 - \rho(A) = d(A) -\iota(A).
\eeq

We will illustrate all these invariants in Example~\ref{ex:complete} below.
  
By abuse of notation we will refer to $\theta(A)$ as the {\em defect} of the configuration $A$ or of a Gale dual $B$. 
As every homogeneous configuration $A$ may be viewed as defining 
a Cayley decomposition with $r=0$ and conditions i) and ii) in Theorem~\ref{th:fi} are satisfied taking $S=\{0\}$, we have that
$\theta(A) \ge 0$ for all $A$. 
If $\theta(A)>0$ we will say that $A$ (or $B$) is {\em defective}.

\begin{remark}\label{rem:3}
Since all the above invariants are defined purely in terms of the configuration $A$ it should be possible to prove the leftmost equality in (\ref{eq:equalities}) 
with purely linear algebraic arguments, as we did for the other equalities in Sections~\ref{sec:left=} and \ref{sec:right=}. 
 However we have only been able to show, with such direct methods, that 
$\theta(A) \leq m(A) - 1 - \lambda(A)$ (cf. Proposition~\ref{ito-flag1}).  We point out, however, 
that the following example shows that~(\ref{eq:equalities}) 
does not necessarily hold for configurations defined over fields of positive characteristic.
\end{remark}

\begin{example}\label{ex:char2}
Let $\k={\mathbb F}_2$ be the field with two elements and consider the following  three-dimensional configuration $A \subset \k^7$:
\beq\label{eq:Achar2}
A \ =\ \left(\begin{array}{ccccccc}1 & 1 & 1 & 0 & 0 & 0 & 0 \\0 & 0 & 0 & 1 & 1 & 1 & 1 \\1 & 0 & 0 & 0 & 0 & 1 & 1 \\0 & 1 & 0 & 0 & 1 & 0 & 1\end{array}\right).
\eeq
A Gale dual $B = \{b_1, \dots, b_7\}$ of $A$ over $\k$
consists of all non-zero vectors in $\k^3$, suitably ordered.
All rank one flats in $B$ are non-splitting but, since ${\rm {char}}(\k) = 2$, all rank two flats are splitting. 
Thus, over $\k$,  $\lambda(A) =1$ and $m(A) - 1 - \lambda(A) = 1$.
Indeed, the associated matroid is known as the Fano matroid and it contains exactly $7$ lines containing $3$ vectors adding-up to $0$. 
Each one of these lines $L$ gives rise to a Cayley decomposition $B= B_0\cup B_1$ where $B_0$ has cardinality $3$ (and $L(B_0)=L$).  Via a linear change of coordinates, we can assume that 
 $B_0 = \{ b_1, b_2,b_3\}$ and $ B_1 =\{b_4, b_5,b_6, b_7\}$, which gives the corresponding Cayley decomposition of $A$ which is visible in~(\ref{eq:Achar2}), where $r=1$, $c=2$ and  $r-c=-1$. Thus, the FI-decomposition is just
$A$ itself, with $r=0$, $c=0$ , $r-c=0 = \theta(A) < m(a) -1 - \lambda(A)$, and the leftmost equality in~(\ref{eq:equalities}) does not hold.

On the other hand, if we consider the same matrix $A$ but over $\Z$  (or over any field with characteristic different from $2$), a Gale dual of $A$ is given by the following vectors in $\Z^3$:
$$\{-e_2-e_3, -e_1-e_3, e_1+e_2+2e_3, -e_1-e_2-e_3, e_1, e_2, e_3\}.$$
It is then clear that $\{e_1\} \subset \{e_1,e_2\}$ defines a non-splitting flag of length $2$. It is still true that
 the FI-decomposition is again given by
$A$ itself, with $r=0$, $c=0$ and $r-c= \theta(A) = m(a) -1 - \lambda(A) = 3-1-2 =0$. Hence, $X_A$ is non-defective over $\C$.
\end{example}

\medskip

\subsection{\bf{The Gale dual setting}} \label{ssec:Gale}

We will now translate these concepts to the Gale dual setting for an arbitrary homogeneous configuration. 

 Let $A \subset  \k^{d+1}$ be a non-pyramidal 
homogeneous configuration such that $d(A) = d$. Suppose $A = A_0 \ast \cdots \ast A_r$ and a $\k$-linear subspace $S$ 
 of $\k^{d+1}$ of dimension $c$ define a decomposition satisfying i) in Theorem~\ref{th:fi} (with $\Z$ replaced by $\k$).
In particular, up to left multiplication by an invertible matrix and permutation of columns, the associated matrix $A$ may be written as 
\begin{equation}
\label{eq:matrixCclower}
 A = \left(\begin{array}{cccc}
A'_0& { 0} & \cdots & { 0}\\
{ 0}  & A'_1 & \cdots & { 0}  \\
{0}  & { 0}  & \cdots & A'_r \\
C_0& C_1 &\cdots & C_r
\end{array}\right),
\end{equation}
where the $A'_i$ are homogeneous configurations in 
$\k^{d_i+1}$  with $d(A'_i)=d_i$ , $\sum_{i=0}^r d_i = d+1-c$,  and the lower matrix $C = (C_0\,  C_1 \dots C_r) \in \k^{c \times n}$  has rank $c$.  
 
 The following result will be needed below:
 
\begin{lemma}\label{lemma:cbarc}
With notation and hypotheses as above, assume moreover that $S$ is minimal in the sense that we cannot find a subspace $S'$ properly contained in $S$ so that  i)  in Theorem~\ref{th:fi} remains valid.
Then, for any index $j, \ j=0,\dots,r$, we have that
\beq\label{eq:cbarc}
\rank\left(\begin{array}{cccccc}C_0 & \cdots & C_{j-1} & C_{j+1} & \cdots & C_r\end{array}\right) = c.
\eeq
\end{lemma}

\begin{proof} 
Assume that $j=0$ and that the matrix $C'= (C_1 \dots C_r)$ has rank $c'<c$.  We may identify $S$ with 
the span of the last $c$ vectors in the chosen basis of $\k^{d+1}$. 
Then, the span $S'$  of the columns of $C'$ has rank $c'<c$, and modding out by $S'$ results in a configuration of join type of the form
\begin{equation}
\label{eq:matrixCclower2}
 \left(\begin{array}{cccc}
A'_0& { 0} & \cdots & { 0}\\
{ 0}  & A'_1 & \cdots & { 0}  \\
{0}  & { 0}  & \cdots & A'_r \\
C'_0& 0 &\cdots & 0
\end{array}\right),
\end{equation}
with $C'_0$ of rank $c-c'>0$. This contradicts the minimality of $S$.
\end{proof}

Given
 a Gale dual of a configuration  $B \subset \k^m$ of $A=A_0 \ast \cdots \ast A_r$, the elements $b \in B$, corresponding to the $A_j$-columns of $A$  
\begin{equation}\label{eq:Bj}
B_j = \{b_i \in B : a_i \in \{e_j\}\times A_j\},
\end{equation}
define a dual-homogeneous subconfiguration $B_j$ and we get a decomposition:
\beq\label{eq:bdecomp}
B = B_0 \cup B_1 \cup \cdots B_r.
\eeq
Conversely, any decomposition of a Gale dual $B$ of $A$  
into dual-homogeneous subconfigurations as in (\ref{eq:bdecomp})
defines a corresponding Cayley decomposition
$A = A_0\ast A_1 \ast \cdots \ast A_r$
of length $r$.  
Moreover, $A$ is equivalent, possible after renumbering, to a configuration of the form (\ref{eq:matrixCclower}), 
where the row span of $C$ gives the affine relations among the elements in $B$ 
which involve elements from at least two distinct $B_j$.

\smallskip

This may be described in more intrinsic terms as follows. Assume $A=A_0 \ast \cdots \ast A_r$ and $S \subset\k^{d+1}$ is a subspace of dimension $c$
which define a decomposition of $A$ satisfying  i) of Theorem~\ref{th:fi} (replacing $\Z$ with $\k$) with $S$ minimal. Let
$\pi: \Z^{d+1} \to \Z^{d+1}/S$. 
 Consider the short exact sequences as in (\ref{eq:exseq}) and the commutative diagram:
\beq\label{fi:gale}
\begin{tikzcd}
0 \arrow[r] & W \arrow[d,hook, "\rho"] \arrow[r,"\iota"] & \k^n \arrow[d,"{\mathrm{id}}"] \arrow[r,"\alpha"] & \k^{d+1} \arrow[d,"\pi"] \arrow[r] & 0\\
0 \arrow[r] & \hat W  \arrow[r,"\hat \iota"] & \k^n  \arrow[r,"\hat\alpha"] & \k^{d+1}/S  \arrow[r] & 0,
\end{tikzcd}
\eeq
where $W = \ker_\k(A)$ and $\hat W = \ker_\k \pi(A)$.  Dualizing, we get 
\beq\label{fi:galedual}
\begin{tikzcd}
0 \arrow[r] & (\k^{d+1}/S)^* \arrow[d, "\pi^*"] \arrow[r,"\hat\alpha^*"] & (\k^n)^* \arrow[d,"{\mathrm{id}^*}"] \arrow[r,"\hat\iota^*"] & \hat W^* \arrow[d,"\rho^*"] \arrow[r] & 0\\
0 \arrow[r] & (\k^{d+1})^*  \arrow[r, "\alpha^*"] & (\k^n)^*  \arrow[r,"\iota^*"] & W^*  \arrow[r] & 0
\end{tikzcd}
\eeq
The map $\rho^*$ is surjective and its kernel is isomorphic to the row span of the matrix $C$ in (\ref{eq:matrixCclower}) via the map:
$$(u_1,\dots,u_n)\in {\mathrm{rowspan}}(C) \mapsto \sum_{i=1}^n u_i\,\hat\iota^*(\xi_i),$$
where as before $\xi_1,\dots,\xi_n$ denotes the dual of the standard basis of $\k^n$.  

Since $\pi(A)$ is a join, a Gale dual $\hat B$ of $\pi(A)$ may be written as 
$\hat B = \hat B_0 \cup \cdots \cup \hat B_r,$
where $\hat B_j$ is a Gale dual of the configuration $A'_j$ depicted in~(\ref{eq:matrixCclower2}) and
$L(\hat B) = L(\hat B_0) \oplus \cdots \oplus L(\hat B_r)$.
In particular, since $\dim L(\hat B) = \dim L(B) + c$, we have
\beq\label{eq:cid1}
c = \sum_{j=0}^r \dim(L(\hat B_j)) - \dim L(B).
\eeq


The following result is now a consequence of  (\ref{eq:cid1}) together with Lemma~\ref{lemma:cbarc}:

\begin{proposition}\label{lem:ABC}  
Suppose $A$ is a non-pyramidal homogeneous Cayley configuration $A = A_0 \ast \cdots \ast A_r$,  $B$ is a Gale dual of $A$,  and $B= B_0 \cup \dots \cup B_r$ is the corresponding Cayley decomposition  as in~(\ref{eq:Bj}). Let $c$ be the maximal
dimension of a subspace $S \subset \k^{d+1}$ such that under the projection to the quotient by $S$ the image $\pi(A) = \pi(A_0)\ast \cdots \ast \pi(A_r)$ is a join.
Then, for each $j=0,\dots,r$, $B_j$ is dual-homogeneous and 
\begin{equation}\label{eq:cdual}
 c = \sum_{j=0}^r \dim L(B_j)  - \dim L(B).
\end{equation}
\end{proposition}

\bp
Because of (\ref{eq:cid1}) we only need to show that $\dim L(B_j) = \dim L(\hat B_j)$ for $j = 0,\dots,r$.  But this follows from the fact that for each $j =0,\dots,r$, $\rho^*(\hat B_j) = B_j$ and the map
$$\rho^* \colon L(\hat B_j) \to L(B)$$
is injective. In fact, take $j=0$, and assume that 
$L(\hat B_0)$ is spanned by $\hat \beta(\xi_1), \dots, \hat \beta(\xi_{n_0})$.  Suppose there exists an element 
$u_1 \,\hat \iota^*(\xi_1) + \cdots + u_{n_0}\,\hat \iota^*(\xi_{n_0}) \in \ker(\rho^*)$, where $n_0 = n(A_0)$.
Then $(u_1,\dots,u_{n_0},0,\dots,0) \in {\mathrm{rowspan}}(C)$, 
implying that $$\rank(C_1,\dots,C_r)<c$$ and thus contradicting Lemma~\ref{lemma:cbarc}.
\end{proof}

The following examples illustrate the constructions above.

\begin{example} Consider the configuration 
$$
A =  \left(\begin{array}{cccccc}
1 & 1 &1 & 0 & 0&0\\
0 & 0&0 &1&1&1\\
1 & 0&0 &1&-1&0\\
\end{array}\right).$$
A Gale dual of $A$ is given by the configuration $B \subset \Q^3$:
$$B = \{b_1,\dots,b_6\}=\{e_1, e_2, -e_1 - e_2, e_3, e_1 + e_3, -e_1 - 2e_3\}.$$
$B$ has only two proper dual-homogeneous subconfigurations: 
$B_0 =  \{b_1,b_2,b_3\}$ and $ B_1 =  \{b_4,b_5,b_6\}$.  Hence the only  possible Cayley decompositions of $B$ are
$B$ and $B_0 \cup B_1$.  Since $c(B_0,B_1) = 1$, both Cayley decompositions are $\theta$-decompositions and $B$ itself is the FI-decomposition.  
Note that the matrix $A$ is in the form (\ref{eq:matrixCclower}) and the bottom row gives the relation:
$$b_1 + b_4 - b_5 =0.$$
\end{example}

\begin{example} \label{ex:return}
We return to Example~\ref{ex:itcirc}.  
For the case of the octahedron $A = \overline{A'}$, where $A'$ is described in (\ref{octahedron}), a Gale dual $B$ is given by 
$$B = \{b_1,\dots,b_6\}=\{e_1, e_1, e_2, e_2, -e_1 - e_2,  -e_1 - e_2\}.$$
The only proper dual-homogeneous subconfigurations correspond to the choice of one of each of $\{e_1,e_2,-e_1-e_2\}$. 
Thus we can obtain decompositions 
$B = B_0 \cup B_1$ with $r=1$ and $c=2$.  The FI-decomposition is then $B$ itself and $B$ is non-defective.

On the other hand in the case of the prism $A = \overline{A'}$ with $A'$ is described in (\ref{prism}), a Gale dual $B$ is given by
the configuration $B \subset \Q^2$:
$$B = \{b_1,\dots,b_6\}=\{e_1, e_2, -e_1 - e_2, -e_1, -e_2, e_1 + e_2\}.$$
There are eight proper dual-homogeneous subconfigurations of $B$:
$$B_0 =  \{b_1,b_2,b_3\}\,;\ B_1 =  \{b_4,b_5,b_6\}\,;\ C_0 =\{b_1,b_4\}\,;\ C_1 =\{b_2,b_5\}\,;\ C_0 =\{b_3,b_6\},$$
together with $D_{ij} = C_i \cup C_j$, for $i\not=j$.
Thus, the possible Cayley decompositions of $B$ are
$$B,\ B_0 \cup B_1,\ C_0 \cup C_1 \cup C_2,\ D_{ij}\cup C_k,$$
for distinct indices $i,j,k$.  It is now easy to see that $\theta(B)=0$ and
$$\theta(B_0,B_1) = 1-2=-1\,;\quad \theta(C_0,C_1,C_2) = 2-1=1\,;\quad \theta(D_{ij},C_k) = 1-1=0.$$
Thus, $B = C_0 \cup C_1 \cup C_2$ is the only $\theta$-decomposition and therefore it  is the FI-decomposition. 
We have $\theta(B) = 1$ and, since all rank-one flats are splitting, $\lambda(B) = 0$. We then have that
$\theta(B) = m(B) -1 - \lambda(B)=1$. 

\end{example}

We next consider Example~5.8 in \cite{FI}:

\begin{example} \label{ex:complete}
The configuration $A=A_0*A_1*A_2*A_3 \subset \Z^{6}$ is the Cayley configuration of length $3$ associated with the lattice configurations in $\Z^2$ 
$$A_0 = \{0, (1,0), (2,0)\}, \, A_1 = \{ 0, (0,1), (0,2)\}, \, A_3 = A_4 = \{ 0, (1,0), (0,1), (1,1)\}.$$
We will use the Gale dual formulation to show that this Cayley decomposition is the FI-decomposition of $A$.
The configuration 
\begin{equation}\label{eq:Bito}
B = B_0 \cup B_1 \cup B_2 \cup B_3,
\end{equation}
where, in terms of the canonical basis $\{e_1,\dots,e_8\}$ of $\Q^8$:
$$B_0 = \{e_1-e_8, -2e_1+e_8,e_1\}\ ;\quad
B_1 = \{e_2-e_7-e_8, -2e_2+e_7+e_8,e_2\}\ ;$$
$$B_2 = \{e_3 +e_5+e_6+e_8, -e_3-e_5, -e_3 - e_6,e_3 -e_8\}\ ;$$
$$
B_3 = \{e_4 -e_5-e_6+e_7, -e_4+e_5, -e_4 + e_6 -e_7, e_4\}
$$
is a Gale dual of $A$ and this decomposition is associated to the given Cayley decomposition of $A$.
We denote the elements of $B$ by $\{b_1,\dots,b_{14}\}$ according to the listing above. 
We claim that the only dual-homogeneous 
 flats in $B$ must be unions of the $B_i$'s. Let $H$ be a dual-homogeneous flat, its projection  to
 $L(f_1)$ must be dual-homogeneous and is contained  
in the projection of  $B$ which equals $\{e_1, -2 e_1, e_1, 0, \dots, 0\}$, where the non-zero terms come from the projection of $B_0$. Thus,  if $H$ contains
 any element in $B_0$ then it must contain all of $B_0$. Similarly, considering the projections to $L(e_2), L(\{e_3,e_5,e_6,e_8\})$ and
 $L(\{e_4,e_5,e_6,e_7\})$ we see that for any $i=1,2,3$, $H \cap B_i \not =\emptyset$ implies $B_i \subset H $.  
 This proves our claim.
  Now, we may verify that for all $i\not=j$, 
$c(B_i,B_j)=0$ and therefore $\theta(B_i \cup B_j) \geq 1$.
Therefore, by Corollary~\ref{fi:cor2}, if
$i\not=j$, $B_i \cup B_j$ may not be a term in an FI-decomposition.   
Similarly we may check that 
$$c(B_0,B_1,B_2) = c(B_0,B_1,B_3) = 0\ ;\quad c(B_0,B_2,B_3) = c(B_1,B_2,B_3)=1.$$
Hence, any union of three distinct $B_i$ is defective and cannot be a term in an FI-decomposition.
Thus, the only dual-homogeneous flats that can be part of an FI-decomposition are the $B_i$'s themselves.
This implies that the Cayley decomposition~(\ref{eq:Bito})
is the (unique)  FI-decomposition. Note that by Proposition~\ref{lem:ABC} we have
$c=\sum_{j=0}^3 \dim L(B_j) - \dim L(B)= 10 -8=2$. Thus, $\theta(A)= 3-2=1$.  From (\ref{eq:equalities}) it follows that 
$\iota(A) = 4$,
$\lambda(A) = 6$,  and $\rho(A) =12$.
 
The subset
 $Z =\{a_1,a_2,a_4,a_5,a_7,a_8\}$
 defines a $4$-dimensional circuit, and thus an iterated circuit of maximal dimension.  
 More generally we make a four-dimensional circuit by 
 choosing any two points from $A_0$, any two points from $A_1$ and any two points from either $A_2$ or $A_3$.  
 
 We can construct a  non-splitting flag of length $6$ as follows:
 Note that $$L(Z) \cap A = Z \cup \{a_3,a_6,a_9,a_{10}\}$$
 and therefore it follows from (i) in Lemma~\ref{circuit-lemma} -as well as by direct computation- that the corresponding elements in $B$ are linearly independent.
 Thus, we can take any non-splitting flag of length $2$ in $B_3$ such as 
 $\{b_{11}\} \subset \{b_{11},b_{12}\} $, 
 and add one by one the linearly independent elements $b_4,b_6,b_9,b_{10}$ to produce a non-splitting flag of length $6$.
 
   For reasons of economy of space we will not write a matrix $M_{\csig}(A)$ with 
 rank equal to $\rho(A)=12$, but will instead observe that we can easily construct the chain 
 $\csig$ from the description of the maximal dimensional circuit $Z$ as in Section~\ref{sec:left=}.
\end{example} 

\medskip

We will now apply the Gale dual formulation of the Furukawa-Ito construction to obtain some basic properties of $\theta$ and FI-decompositions.  
We begin by showing that these properties are inherited by 
sub-decompositions.  More precisely,

\begin{proposition}\label{fi:lemma3}
Let $B =  B_0 \cup \cdots \cup B_r$ be a $\theta$-decomposition (resp. an FI-decomposition).  Then for every subset 
$I = \{i_0,\dots,i_\ell\}  \subset \{0,\dots,r\}$,
the decomposition
$$B_I = B_{i_0}\cup \cdots \cup B_{i_\ell}$$
is a $\theta$ (resp. FI)-decomposition.
\end{proposition}

\begin{proof}
Assume that   $B =  B_0 \cup \cdots \cup B_r$ is a $\theta$-decomposition, $m = \dim(L(B))$, and that $I = \{0,\dots,\ell\}$.
If $B_I =  B_0 \cup \cdots \cup B_\ell$ is not a $\theta$-decomposition then there exists another Cayley decomposition: 
$B_I = C_0\cup \cdots \cup C_k$  such that
$$k - c(C_0,\dots,C_k) > \ell - c(B_0,\dots,B_\ell).$$ 
Consider the Cayley decomposition of length $r' = k + r - \ell$: $$B = C_0 \cup \dots \cup C_k \cup B_{\ell+1} \cup \dots \cup B_r,$$
 which we
will just write as $B=C_0 \cup \dots \cup B_r$ and set 
$c' = c(C_0,\dots,B_r)$. Then, 
\beas
r' - c' &=& r' - \sum_{j=0}^k \dim(L(C_j)) - \sum_{j=\ell+1}^{r} \dim(L(B_j)) + m\\
&=&  (r - \ell) + (k - c(C_0,\dots,C_k) +\dim(L(B_I))  - \sum_{j=\ell+1}^{r} \dim(L(B_j)) + m\\
&>& (r - \ell) + (\ell - c(B_0,\dots,B_\ell) +\dim(L(B_I))  - \sum_{j=\ell+1}^{r} \dim(L(B_j)) + m\\
&=& r - \sum_{j=0}^{\ell} \dim(L(B_j)) - \sum_{j=\ell+1}^{r} \dim(L(B_j)) + m
= r - c(B_0,\dots,B_r), 
\eeas
contradicting the assumption that $B =  B_0 \cup \cdots \cup B_r$ is a $\theta$-decomposition.

Finally, if we assume that   $B =  B_0 \cup \cdots \cup B_r$ is an FI-decomposition, then we already know that 
$B_I = B_0 \cup \dots \cup B_\ell$ is
a $\theta$-decomposition. If $B_I$ admitted a $\theta$-decomposition of length $k<\ell$, then this would yield a new  decomposition which, by an argument
 completely analogous to the one above, is seen to be a $\theta$-decomposition. But this $\theta$-decomposition  would have fewer than $r$ terms which is impossible.

\end{proof}

Applying Proposition~\ref{fi:lemma3} to the case $I=\{j\}$ we get:

\begin{corollary}\label{fi:cor2} 
If $B = B_0 \cup \cdots \cup B_r$ is a $\theta$-decomposition, then 
for every $j=0,\dots,r$, $\theta(B_j) = 0$.
\end{corollary}

Moreover we may also deduce the following

\begin{corollary}\label{fi:cor3} 
Suppose $B = B_0 \cup \cdots \cup B_r$ is a $\theta$-decomposition.  Let 
$I = \{i_0,\dots,i_\ell\}  \subset \{0,\dots,r\}$, $|I| \geq 2$, and set
$B_I = \bigcup_{i\in I} B_i$.  Then
\beq\label{fi:herineq}
\sum_{i\in I} \dim(L(B_i)) - \dim(L(B_I)) \leq |I| - 1.
\eeq
Moreover, if $B = B_0 \cup \cdots \cup B_r$ is the FI-decomposition then the inequality above is strict and, in this case
\begin{enumerate}
\item[i)] For each $0\leq i,j \leq r$, $i\not= j$, $L(B_i) \cap L(B_j)=0$.
\item[ii)] For every $j=0,\dots,r$, $B_j \subset B$ is a flat.
\end{enumerate}
\end{corollary}

\begin{proof} 
The first statement follows from the fact that if $B = B_0 \cup \cdots \cup B_r$ is a $\theta$-decomposition then so is $B_I = \bigcup_{i\in I} B_i$ and, therefore, 
this decomposition computes $\theta(B_I) \geq 0$ which yields inequality (\ref{fi:herineq}).  If $B = B_0 \cup \cdots \cup B_r$ 
is the FI-decomposition then $\theta(B_I) > 0$ and this gives the strict inequality.

Item i) follows from taking $I = \{i,j\}$.  Finally, suppose $B_i$ is not a flat for some $i=0, \dots, r$.  Then, there exist an element $b \in B$ such that 
$b \in  L(B_i) \setminus B_i$. 
But then $b \in B_j$ with $j\not=i$ implying that
$ L(B_i) \cap L(B_j) \not=\{0\}$,
contradicting the previous statement.
\end{proof}

\begin{remark}\label{rem:4}
Note that  Corollary~\ref{fi:cor3}
 gives necessary conditions for a Cayley decomposition to be a $\theta$-decomposition or the FI-decomposition.  
 We expect those conditions 
to be sufficient as well.  We also point out that it follows from the uniqueness statement in Theorem~\ref{th:fi} that the FI-decomposition is unique up to reordering and that the passage from a $\theta$-decomposition to the FI-decompositon is accomplished by collecting together sub-decompositions $B_I$ with $\theta(B_I)=0$. \end{remark}

\begin{lemma}\label{ito:lemma_quotient}
Let $B$ be a configuration and $F \subset B$ a non-splitting rank-one flat.  
 Let $\pi \colon  L(B) \to L(B) /L(F) $
be the natural projection and $\tilde B = \pi(B\setminus F)$.  Then 
$$\theta(B) \leq \theta(\tilde B).$$
\end{lemma}

\begin{proof}
Let $B = B_0 \cup \cdots \cup B_r$ be the FI-decomposition of $B$.  Then by Corollary~\ref{fi:cor3}
we may assume  that $F \subset C_0$ and $L(F) \cap \langle C_j \rangle = \{0\}$ for $j\geq 1$.  Note that
$F\neq C_0$ because $\fs(F)\neq 0$.  Then, 
$\tilde B = \pi(B_0 \setminus F) \cup \pi(B_1) \cup \cdots \cup \pi(B_r)$
is a Cayley decomposition  of $\tilde B$ with the same invariants as the decomposition of $B$.  

Hence,
$\theta(\tilde B) \geq r - c(B_0,\dots,B_r) = \theta(B)$.
\end{proof}  

We may now use  Lemma~\ref{ito:lemma_quotient} to show the following inequality.

\begin{proposition}\label{ito-flag1}
Let $B$ be a dual-homogeneous configuration (as always, we assume $0 \not\in B$) over any field $\k$ of characteristic zero.  Then
\beq\label{ineq}
\theta(B) \leq \dim(L(B)) - 1 - \lambda(B).
\eeq
\end{proposition}

\begin{proof}
If $\dim L(B)  =1$ then both sides of the inequality equal zero.  Suppose now that the result holds for
 $\dim L(B) =m-1$ and let $B$ be a configuration of rank $m$.

Let $\lambda = \lambda(B)$ and $F_1 \subset F_\lambda$ a maximal non-splitting flag in $B$. 
 Let $\pi \colon L(B) \to L(B) /  L(F_1)$ be the natural projection and 
 $\tilde B = \pi(B\setminus F_1)$ the projected configuration.  From Lemma~\ref{quotient-flag} we have that 
$\lambda(\tilde B) = \lambda(B) - 1$.

Now, by inductive hypothesis we may assume that 
$\theta(\tilde B) \leq \dim(L(\tilde B)) - 1 - \lambda(\tilde B)$.
But $\dim(L(\tilde B))= \dim(L(B)) -1$ and by Lemma~\ref{ito:lemma_quotient} we have 
$\theta(\tilde B) \geq \theta(B)$.
Hence, $\theta(B) \leq \dim(L(B)) - 1 - \lambda(B)$, 
as asserted.
\end{proof}

\medskip

We next show how to compute the invariant $\lambda$  (the maximal length of a non-splitting flag)  from a $\theta$-decomposition of $B$.

\begin{theorem}\label{th:fi_dual}
Let $A \subset \Z^{d+1}$ be a homogeneous, non-pyramidal configuration of affine dimension $d$ and $B$ a Gale dual of $A$. 
If 
$  B = B_0\cup B_1 \cup \cdots \cup B_r
$
is any $\theta$-decomposition, we have
\beq\label{eq:max_lambda}
\lambda(A) =  \sum_{j=0}^r (\dim L(B_j) - 1).
\eeq
\end{theorem}

\bp
We have by the leftmost equality in~(\ref{eq:equalities}) that $\lambda = \lambda(A) = m - 1 - \theta(A)$, with $m=m(A) = \dim(L(B))$.
Then, by Proposition~\ref{lem:ABC}, we get
$$ \lambda = m -1 -r +c =  -1 - r + \sum_{j=0}^r (\dim L(B_j)) = \sum_{j=0}^r (\dim L(B_j) - 1),$$
as claimed.
\ep

\subsection{\bf{Consequences of 
Theorem~\ref{th:fi_dual}}}\label{ssec:BFI}
 
Recall from \cite[Definition~5]{CC} that a configuration $B\subset \Q^m$ is called {\em irreducible} if any two elements in $B$ are linearly independent.  
If $B$ is non-pyramidal and irreducible then the linear matroid $B$ is {\em simple}. We carry this terminology to a configuration $A \subset \Z^{d+1}$ and 
say that it is irreducible if and only if any (every) Gale dual of $A$ is irreducible.  Recall that a rank-one flat $F \subset B$ is called {\em splitting} if 
${\mathfrak s}(F) = 0$ and {\em non-splitting} otherwise.  Given $B\subset \Q^d$, we denote by $B^{{\rm red}}$ the irreducible configuration 
obtained by eliminating all elements lying in splitting lines and replacing those in non-splitting lines $F$ by
${\mathfrak s}(F)$.  We call $B^{{\rm red}}$  the reduction of $B$.  It is easy to see that $\lambda(B) = \lambda(B^{{\rm red}})$ and therefore 
it follows from (\ref{eq:equalities}) that:
\beq\label{eq:def_red}
\theta(B) = \theta(B^{{\rm red}}) + (\dim L(B) - \dim L(B^{{\rm red}})).
\eeq
Thus, in order to understand the dual defect of toric varieties, it suffices to consider irreducible configurations.  

\begin{remark} We point out that even though, for simplicity, we have deduced (\ref{eq:def_red}) from (\ref{eq:equalities}), 
it can also be deduced directly with linear-algebraic arguments.  
\end{remark}

The following result bounds the dual defect of irreducible configurations.

\begin{proposition}\label{prop:max_def} 
Let $A \subset \Z^{d+1}$ be a homogeneous irreducible configuration with $d(A)=d$.  Then $\theta(A) \leq (m(A) - 2)/2$.
Moreover for the cases of maximal possible defect we have:
\begin{itemize} 
\item
If $m(A) = 2k$ is even and $\theta(A) =k-1$, then $A$ is affinely equivalent to a join of $k$ non-defective configurations $A_0,\dots,A_{k-1}$ with $m(A_j) = 2$.
\item If $m(A) = 2k+1$ is odd and $\theta(A) =k-1$, then either $A$ is affinely equivalent to a join of $k$ non-defective configurations $A_0,\dots,A_{k-1}$ 
with $m(A_0)=3$ and $m(A_j) = 2$, $j\geq 1$, or to a Cayley configuration 
$ A_0 \ast \cdots \ast A_k$
of $k+1$ configurations of codimension $2$ and with $c(A_0,\dots,A_k) = 1$.
\end{itemize}
\end{proposition}

\bp
Let  $B \subset \Q^m$ be a Gale dual of $A$ and 
$B = B_0 \cup \cdots \cup B_r$
the FI-decomposition.  Setting $\lambda_i = \dim L(B_i) -1$, we have by Theorem~\ref{th:fi_dual} that 
\beq\label{eq:partition}
\lambda = \lambda_0 + \cdots + \lambda_r.
\eeq
Moreover, the irreducibility of $A$ implies that $\lambda_i \geq 1$ for all $i$.

Now, since $\theta(A) = r - c \leq r$ we have
$$\theta(A) + 1 \leq r + 1 \leq \sum_{j=0}^r \lambda_j = \lambda = m - 1 - \theta(A).$$
Thus $2 \theta(A) \leq m - 2$ and the result follows.

In particular we have, if $\theta>0$ that
\beq\label{eq:bound_r}1\leq r \leq \lambda - 1 = m - \theta - 2,
\eeq
where the left inequality 
follows from  Corollary~\ref{fi:cor2} that $r\geq 1$.

Suppose now that $m = 2k$ and $\theta(A) = k-1$.  Then $\lambda = m - 1 - \theta = k$ and therefore $r \leq k-1$.  Hence  $\theta(A) = k-1$ 
if and only if $r=k-1$, which implies $\lambda_i =1$ for all $i=0,\dots,r-1$, and $c=0$.  Hence $A$ is equivalent to a join as asserted.

Finally, suppose $m = 2k+1 $ and $\theta(A) = k-1$.  Then 
$\lambda = 2k - k + 1 = k+1$ and there are two possible partitions of $\lambda$ of length at least $k$:
$$\lambda_0=\cdots=\lambda_k =1\ ;\hbox{\ and\ }
\lambda_0=2\ , \lambda_1=\cdots=\lambda_{k-1} =1.$$
In the first case we must have $c=1$ and in the second $c=0$.  This yields the result.
\ep

The following consequences of Proposition~\ref{prop:max_def} are expressed in terms of Gale duality to facilitate the comparison with results in \cite{CC}

\begin{example}\label{m=3}
If $m=3$, then $\theta < 1$ and therefore, every irreducible configuration with $m(A) = 3$ is non-defective. This is Theorem~20 in \cite{CC}.  
\end{example}

\begin{example}\label{m=4}
If $m=4$ then $\theta \leq 1$ and, therefore by Proposition~\ref{prop:max_def}, if $B$ is defective and irreducible then  $B = B_0 \cup B_1$, 
where $\dim L(B_0) = \dim L(B_1)$ and $L(B_0) \cap L(B_1) = \{0\}$.  This is Theorem~21 in \cite{CC}.
\end{example}

\begin{example}\label{m=5}
If $m=5$ then, again $\theta \leq 1$ but now there are two non-equivalent irreducible defective configurations:
\begin{itemize}
\item $B = B_0 \cup B_1$, $\dim L(B_0) = 3$, $\dim L(B_1) = 2$.
\item $B = B_0 \cup B_1 \cup B_2$,  $\dim L(B_i) = 2$, for all $i$ and $L(B_i) \cap L(B_j) = \{0\}$ for $i\not =j$.
\end{itemize}
We point out that, although possible, obtaining this result with the methods in \cite{CC} would involve very long and awkward arguments.
\end{example}

It follows from Proposition~\ref{prop:max_def} that if  $\theta(A) > (m(A) - 2)/2$ then $A$ may not be irreducible.  We can give a more precise statement.

\begin{proposition}\label{prop:excess}
Let $A \subset \Z^{d+1}$ be a homogeneous configuration of maximal rank and suppose that $\theta(A) > (m(A) - 2)/2$.  Then, 
\begin{itemize} 
\item
If $m(A) = 2k$ and $\theta(A) =k-1 + \ell$, then 
$$\dim L(B) - \dim L(B^{{\rm red}}) \geq 2\ell.$$
\item If $m(A) = 2k+1$  and $\theta(A) =k-1 + \ell$, then 
$$\dim L(B) - \dim L(B^{{\rm red}}) \geq 2\ell -1.$$
\end{itemize}

\end{proposition}

\bp
We prove the even case and leave the analogous odd case to the reader.  Suppose $m=2k$ and $\theta(B) = k - 1 + \ell$.  Let $m_1 = \dim L(B^{{\rm red}})$, then
$$ k - 1 + \ell = \theta(B) = \theta(B^{{\rm red}}) + (m - m_1)$$
$$(m - m_1) = \theta(B) - \theta(B^{{\rm red}}) \geq k - 1 + \ell - (m_1 -2)/2 = \ell + (m - m_1)/2$$
and therefore
$m - m_1 \geq 2\ell$ as asserted.
\ep

Given a configuration $A$ we have $\theta(A) \leq m(A)-1$.  If we have equality, then $\lambda(A) =0$ and $\dim X_A = \dim X^{\vee}_A$ 
and it follows that  every rank-one flat in $B$ must be splitting or, equivalently, that $\dim L(B^{{\rm red}})=0$ (cf. \cite{BDR}).  Moreover, as shown in 
\cite[Theorem~3.3]{BDR}, if $X_A$ is equivariantly embedded then $X_A$ is self-dual; i.e. $X_A \cong X^{\vee}_A$.

\begin{proposition}\label{prop:big_defect}
Let $A$ be a configuration with $\theta(A) = m(A) -2$ and let $B$ be a Gale dual.  Then $B^{{\rm red}}$ is a non-defective configuration of rank $2$.  
Similarly, if $\theta(A) = m(A) -3$, then either $B^{{\rm red}}$ is a non-defective configuration of rank $3$ or it is a rank $4$ configuration with defect equal to $1$.
\end{proposition}

\bp
In the first case we must have $\lambda = 1$ therefore the FI-decomposition of $B$ must be of the form
$B = \Lambda_1 \cup \cdots \cup \Lambda_s \cup B_0$,
where $\Lambda_i$ are dual-homogeneous rank-one flats and $m(B_0) = 2$, 
$\lambda(B_0)=1$.  Clearly $B^{{\rm red}} = B_0^{{\rm red}}$.

If $\theta(A) = m(A) -3$  then $\lambda(A)=2$ and therefore there are two possibilities for the FI-decomposition of $B$:
$$B = \Lambda_1 \cup \cdots \cup \Lambda_s \cup B_0, \, \hbox{ \ or } \,
 B = \Lambda_1 \cup \cdots \cup \Lambda_s \cup B_1\cup B_2,$$
 with $\dim L(\Lambda_i) = 1$, $\dim L(B_0) = 3$ and $\dim L(B_1) = \dim L(B_2)=2$.  Since 
in the first case $B^{{\rm red}} = B_0^{{\rm red}}$ and in the second $B^{{\rm red}} = (B_1\cup B_2)^{{\rm red}}$, the result follows.
\ep

\end{document}